\RequirePackage{fix-cm}
\documentclass[smallextended]{svjour3}       % onecolumn (second format)
\usepackage[margin=1in]{geometry}
\smartqed  % flush right qed marks, e.g. at end of proof

\usepackage{cite}
\usepackage{graphicx}
\usepackage{amsmath,amssymb,amsfonts,latexsym}
\usepackage{epsfig}
\usepackage{subfigure}
\usepackage{algorithm}
\usepackage{array}
\usepackage{url}
\usepackage{float}
\usepackage{enumitem}
\usepackage{caption}
\usepackage{algpseudocode}

\usepackage{color}
\newcommand{\bl}[1]{\textcolor{black}{#1}}

\def \vx{ \mathbf{x} }

\newcommand{\lf}{\left}
\newcommand{\rt}{\right}

\def \vx{ \boldsymbol{x} }

\def \vy{ \boldsymbol{y} }

\def \vx{ \boldsymbol{x} }

\journalname{JSC}
\begin{document}

\title{Curvilinear Mesh Adaptation using Radial Basis Function Interpolation and Smoothing}

\titlerunning{Curvilinear Mesh Adaptation using RBF Interpolation and Smoothing}        % if too long for running head

\author{Vidhi Zala         \and
        Varun Shankar  \and
        Shankar P. Sastry \and
        Robert M. Kirby %etc.
}

%\authorrunning{Short form of author list} % if too long for running head

\institute{Vidhi Zala \at
              Scientific Computing and Imaging Institute, University of Utah, Salt Lake City, UT 
							\email{vidhi@sci.utah.edu}
%             \emph{Present address:} of F. Author  %  if needed
           \and
           Varun Shankar \at
              Department of Mathematics, University of Utah, Salt Lake City, UT\\
							\email{vshankar@math.utah.edu}
              \and
              Shankar P. Sastry \at
                Scientific Computing and Imaging Institute, University of Utah, Salt Lake City, UT
              \and
              Robert M. Kirby \at
               Scientific Computing and Imaging Institute, University of Utah, Salt Lake City, UT   \\     
               \email{kirby@sci.utah.edu}
}

\date{Received: date / Accepted: date}
% The correct dates will be entered by the editor

\maketitle
\parindent=0pt
\begin{abstract}
We present a new iterative technique based on radial basis function (RBF) interpolation and smoothing for the generation and smoothing of curvilinear meshes from straight-sided or other curvilinear meshes. Our technique approximates the coordinate deformation maps in both the interior and boundary of the curvilinear output mesh by using only scattered nodes on the boundary of the input mesh as data sites in an interpolation problem. Our technique produces high-quality meshes in the deformed domain even when the deformation maps are singular due to a new iterative algorithm based on modification of the RBF shape parameter. Due to the use of RBF interpolation, our technique is applicable to both 2D and 3D curvilinear mesh generation without significant modification.
\end{abstract}

\section{Introduction}
\label{intro}

The increasing use of simulations built upon high-order numerical methods in practical engineering problems
necessitates the generation of meshes that conform to irregular domain geometries.  To maintain the high-order
numerical nature of these simulations, the geometric accuracy of the domain must also be high-order, thus 
motivating (high-order) curvilinear meshes.  The starting point for many high-order meshing techniques is to create
a valid low-order (straight-sided) mesh which is then ``adapted" to the curved geometry.  The challenge when
accomplishing this update is the balancing act between faithfully representing the boundaries of interest while
maintaining a mesh whose elements of are of good quality (and hence have favorable numerical properties).  We 
present a technique based on RBF-interpolation that produces superior quality meshes by first deforming the 
domain to meet the geometric constraints of the problem and then iteratively adapting or smoothing the
mesh in a way to capitalize on the properties of the RBF-interpolation functions we employ.
In this work, mesh adaptation has been studied in \bl{the} context of both refinement (coarsening) and smoothing: 
the technique presented herein aims to combine the best of both approaches.  

%The fundamental difference between these two techniques is that smoothing preserves the topology of the mesh while attempting to improve the quality by node re-positioning, but refinement/coarsening changes the number of elements inside the mesh by inserting or collapsing existing elements. While our technique aims to progressively improve the element qualities in mesh, it also varies the number of elements by inserting new elements or collapsing the bad ones. Ultimately, the resulting mesh quality, as calculated to be the average of all elemental qualities, is higher than the initial mesh. According to popular definitions of optimization of
%the mesh quality, the refinement and the coarsening mainly
%try to optimize the mesh density, while edge swapping and
%mesh smoothing mainly aim to optimize the shape regularity\cite{imrChen04}. 
%A form of mesh smoothing called optimization-based smoothing is gaining attention in the meshing community. \cite{FreitagDiachin2006}. It can be viewed as a variation of smoothing methods used classically such as laplacian smoothing with the change in objective function to optimize. In case of heuristic based smoothing procedures like laplace smoothing, movement of the mesh nodes are aimed to minimize distortion by placing vertices at some optimal distance from its neighbors. Our technique optimizes an objective function based on the mesh quality metrics in a localized stencil which has global implications. 

We begin by summarizing the state-of-the-art in high-order (curvilinear) mesh generation.
Over the last decade, many techniques for the generation and deformation of linear meshes into curvilinear meshes have been proposed \cite{Moxey2015636,GPRPS13,GPRPS14,Toulorge:2013:RUC:2527809.2527935,POPersson}.  Sastry et al.\cite{Sastry2015135} provided the following taxonomy for partitioning the literature landscape: optimization-based methods \cite{Remacle2013,Sastry2012}, PDE-based meshods \cite{POPersson}, and interpolation-based
techniques \cite{deBoer:2007:MDB:1243516.1243757}. 
Some of the notable work done in the mesh deformation and curvilinear mesh generation 
can be attributed to the application of one or more techniques from these three classes. 
The optimization-based techniques aim at optimizing an objective function depending on 
the geometry of the domain and the mesh.  Sastry et al. \cite{Sastry2012} proposed a
log-barrier optimization routine to dictate vertex movement and to improve the quality of a tangled mesh (due to the
deformation) to obtain a valid mesh. The Remacle group \cite{Remacle2013, Geuzaine2015} 
developed a log-barrier technique that generates a valid mesh by maximizing the minimum Jacobian of high-order elements
in the mesh.  From the class of PDE-based methods, Moxey et al. \cite{Moxey2016130} presented a technique based on the thermo-elastic analogy by modelling the mesh as non linear elastic material.  In a subsequent paper, Moxey et al. \cite{Turner2016340} take the variational approach further by optimizing mesh quality using a scaled Jacobian approach. Sastry et al.\cite{Sastry2015135} compared and contrasted the thermo-elastic 
method with the RBF interpolation using thin-plate-splines. Experiments for that effort help 
establish the superiority of the RBF interpolation-based technique by generating elements that are of higher 
quality and \bl{conform to} the boundary geometry.

The last class, interpolation-based methods, has been mostly applied to linear mesh deformation functions. Staten et al. \cite{Staten2012} developed the simplex-linear transformation algorithm, which carries out a linear interpolation of mesh vertices after making a coarse mesh as initial step. Sastry et al. \cite{Sastry2015135} developed a technique for curvlinear mesh generation using thin-plate spline RBFs, which belong to the class of polyharmonic splines. They further demonstrated that interpolants based on polyharmonic spline help preserve the shape of elements after deformation. However, the thin-plate spline technique did not possess the ability to deal with degenerate deformation maps, or smooth any resulting mesh tangles. Further, the technique did not generalize to 3D meshes in a straightforward fashion.

Broadly speaking, we can treat PDE-based methods and interpolation-based methods as being in the same class, where
the positions of the interior mesh vertices are interpolated from the positions of the boundary 
vertices using either the solution of a PDE or an explicit interpolation technique.  Such a characterization is useful
as it helps motivate our work:  we seek to develop an interpolation-based method that through our choice of 
the interpolating functions mimics some of the favorable properties observed in the PDE-based approaches while being applicable to both 2D and 3D mesh generation.  
We present a generalization of \cite{Sastry2015135} that uses RBFs with a shape parameter to smooth node clusters resulting from singular or non-smooth deformation maps. Specifically, we turn to the Mat\'ern kernels (also referred to as \emph{Sobolov splines}), a family of RBFs closely related to the polyharmonic splines. As their alternate name implies, interpolants based on these kernels are the minimum Sobolev norm interpolants, possessing similar properties to polyharmonic splines, but possessing a shape parameter that is extremely useful for tuning.   
In Section \ref{sec:methods}, we compare the Mat\'ern kernels to the polyharmonic splines, and present a tuning algorithm for the shape parameter to help achieve quasi-local smoothing of these interpolants.

The remainder of the paper is organized as follows. 
In Section \ref{sec:review} we review RBF interpolation with a focus on Mat\'ern kernels, our basis of choice; we justify the use of this basis, and we also present a generalization of existing techniques to smooth RBF interpolants. We go on to present
a mathematical description of our quality heuristics and a new adaptation and smoothing algorithm in 
Section \ref{sec:methods}. We then undertake a thorough complexity analysis of our method in Section \ref{sec:complex}. Finally,
we present numerical experiments exploring the behavior of our method on different classes of deformation functions
in Section \ref{sec:experiments}. We conclude with a discussion of the results and provide some comments on future work.

\section{Review}
\label{sec:review}

\subsection{RBF Interpolation} 
\label{subsec:rbfI}
RBFs are a popular tool for scattered data interpolation in arbitrary dimensions and have become increasingly popular in machine learning \cite{5707295, doi:10.1142/S0129065709002026}, computer graphics \cite{Carr:2001:RRO:383259.383266, 5395267}, mesh generation and repair~\cite{PiretRepair,Sastry2015135} and in the numerical solution of PDEs~\cite{Larsson2003891, SWFKIJNMF2014}. More relevant to this article, RBFs have also been used to interpolate data on co-dimension one submanifolds of $\mathbb{R}^s$ with excellent approximation properties using only straight-line (i.e. Euclidean) distances in the embedding space~\cite{FuselierWright2010}, a feature that has been leveraged to solve PDEs on surfaces~\cite{FuselierWright2012,SWFKJSC2014}. In our application, the relevant submanifolds are the boundaries of (irregular) domains in $\mathbb{R}^2$ and $\mathbb{R}^3$.

We now briefly describe RBF interpolation in $\mathbb{R}^s$; for interpolation on submanifolds $\mathbb{M} \subset \mathbb{R}^s$, it is only necessary for the points to lie on $\mathbb{M}$. Given a set of (scattered) nodes $X = \{\vx_i\}_{i=1}^N$ in $\mathbb{R}^s$ and a set of data values $Y = \{\vy_i\}_{i=1}^N$ sampled from some function $f : \mathbb{R}^s \to \mathbb{R}$, the RBF approximation to $f$ is obtained by a linear combination of \emph{shifts} of a single \emph{radial} kernel or basis function $\phi$ such that
\begin{align}
 \label{eq:rbfI}
I_{\phi}f (\vx,\epsilon) = \sum_{i = 1}^{N} {\lambda}_i(\epsilon) \phi(\epsilon,r_i(\vx))
\end{align}
where $\phi(\epsilon,r_i(\vx)) = \phi(\epsilon\|\vx - \vx_i\|)$ and $\epsilon>0$ is a \emph{shape parameter} that controls the flatness of the RBF. To find the unknown coefficients $\lambda_i$, we enforce the interpolation conditions
\begin{align}
\left. I_{\phi} f \right|_X &= Y, \\
\implies I_{\phi} f (\vx_i,\epsilon) &= \{\vy_i\}_{i=1}^N.
\end{align}
If $\phi$ is a positive-definite radial kernel or an order one conditionally positive-definite kernel on $\mathbb{R}^s$ and all nodes in $X$ are distinct, the above interpolation problem has a unique solution, and the corresponding RBF interpolation matrix is invertible~\cite{Fasshauer:2007}. In the limit as $\epsilon \to 0$ (i.e. a flat kernel), RBF interpolants to data scattered in $\mathbb{R}^s$ typically converge to (multivariate) polynomial interpolants~\cite{DriscollFornberg2002,LarFor05,Sch05}, and to spherical harmonic interpolants on a sphere~\cite{FoPi07}.  For smooth target functions, smaller (but non-zero) values of $\epsilon$ generally lead to more accurate RBF interpolants~\cite{FoWr,LarFor05}. Unfortunately, computing these interpolants by solving the linear system involving the RBF interpolation matrix becomes ill-conditioned for small $\epsilon$ (see, \emph{e.g.,}~\cite{FZ07}).  While some stable algorithms have been developed for bypassing this ill-conditioning~\cite{FoWr,FoPi07,FaMC12,FLF,FoLePo13,FlyerPHS}, these algorithms do not apply when the nodes lie on a lower-dimensional surface than the embedding space. Our approach will be to pick a value of $\epsilon$ that results in some target condition number $\kappa$ in the interpolation matrix that is very close to the edge of ill-conditioning. This typically results in excellent approximation~\cite{SWFKJSC2014}. \bl{Our goal will be to approximate vector-valued functions in this work. We accomplish this by interpolating each component of the vector-valued functions using a scalar RBF interpolant}.

For similar reasons to~\cite{Sastry2015135}, we choose an RBF $\phi$ with global support. \bl{Specifically, we use the piecewise-smooth $C^4$ Mat\'ern kernel given by:
\begin{align}
\label{eq:C4matKer}
\phi(\epsilon r) = (3+3\epsilon r + \epsilon^2 r^2)e^{-\epsilon r}.	
\end{align}}
Our reasons for using this kernel are twofold: first, the reproducing kernel Hilbert space corresponding to this kernel is a standard Sobolev space and therefore well-understood; second, unlike the polyharmonic splines, the Mat\'ern kernel comes equipped with a shape parameter $\epsilon$, such that the limit $\epsilon \to 0$ recovers the polyharmonic spline kernels used in~\cite{Sastry2015135}. Modification of this shape parameter upon evaluation of the RBF interpolant can allow \emph{smoothing}. This will be explained in the following section. For more on Mat\'ern kernels, we refer the reader to~\cite{Fasshauer:2007,FasshauerGreen}.
%\begin{align}
%\label{eq:matKer}
%\phi(\epsilon r) = \frac{K_{\beta-\frac{s}{2}}(\epsilon r) \cdot (\epsilon r)^{\beta-{\frac{s}{2}}}}{2^{\beta-1}\Gamma(\beta)}, \textrm{   }\beta > \frac{s}{2}.
%\end{align}
%Here $K_{\nu}$ represents the \textit{modified Bessel function} of the second kind of order $\frac{s}{2}$. 
%\bl{For all the experiments in Section \ref{sec:experiments}, we use the $C^4$ mat\'ern kernel defined as

%Some simple representations of Mat\'ern kernels are listed in Table \ref{tab:matKer} below for different $\beta$ values.
%% For tables use
%\begin{table}
%% table caption is above the table
%\caption{Mat\'ern kernels for different $\beta$ values.}
%\label{tab:matKer}      % Give a unique label
%% For LaTeX tables use
%\begin{tabular}{ll}
%\hline\noalign{\smallskip}
%%first & second & third  \\
%%\noalign{\smallskip}\hline\noalign{\smallskip}
%{$\beta=\frac{s+1}{2}$} & $e^{-\epsilon r}$\\ 
%{$\beta=\frac{s+3}{2}$} &  $(1+\epsilon r)e^{-\epsilon r}$ \\ 
%{$\beta=\frac{s+5}{2}$} & $-(3+3\epsilon r + \epsilon^2 r^2)e^{-\epsilon r}$  \\
%{$\beta=\frac{s+7}{2}$} & $(15+ 15 \epsilon r + 6 \epsilon^2 r^2 + \epsilon^3 r^3)e^{-\epsilon r} $  \\
%\noalign{\smallskip}\hline
%\end{tabular}
%\end{table}
%
%We note that changing $\beta$ specifies Mat\'ern kernels with different degrees of smoothness. For instance, $\beta = \frac{s+5}{2}$ yields $\phi(r) \in C^4 (\mathbb{R}^s)$. 

\subsection{Mesh Quality}
\label{subsec:MeshQ}

%The vertex re-positioning problem,
%from an optimization perspective leads to
%formulate a single objective function measuring global
%mesh quality. We define the global mesh quality as an aggregate of the 
%qualities of all the elements inside the mesh. The objective function to optimize is
%constructed by accumulating contributions from
%each local stencil quality. A stencil, as used in our technique, is a group of elements surrounding a node in mesh which contributes to the quality heuristic which we use to re-position the node. We
%consider two-fold approach for numerically optimizing
%the global objective function or the overall quality of the mesh: first, evaluating element quality for each element and assigning an aggregate heuristic to nodes based on the neighbor-ring and second, using the heuristic to modify the shape parameter for the RBF-interpolation to recover the deformation map in the interior, thus obtaining a re-positioning of interior nodes which produces better heuristic in the subsequent iterations. The second step in this approach has a global effect, although the shape parameter changes is determined by local stencil. This leads to superior numerical results compared to other techniques that rely on strictly localized or strictly global heuristic optimization. Thus, our contribution to global change in mesh quality based on local change in shape parameters of interpolation is useful and unique.

Our RBF-based technique accomplishes two distinct purposes: first, it recovers deformation maps (and therefore a deformed mesh) using data only on the boundary of the input and output domains; second, it also attempts to automatically smooth the recovered deformation map so as to obtain a deformed mesh with good-quality elements. An element quality metric is a scalar function of node positions
that measures some geometric property of the element~\cite{AMQM2001PKnupp}. In this section, we present a brief overview of the popular metrics for measuring mesh quality. Assume for the following discussion that a mesh contains \bl{a} finite set of vertices V defined as $V = \{\vx_i\}_{i=1}^N$ in $\mathbb{R}^{\bl{s}}$, and a finite set of elements E defined by groupings of those vertices. The elements are triangles in 2D and tetrahedra in 3D.

There are many popular techniques for generating meshes out of point sets, like octree mesh generation~\cite{de1990explicit}, Delaunay triangulation~\cite{george1998delaunay,field1988laplacian,imrChen04} and advancing-front \cite{NME:NME1620382102}. Out of these techniques, the Delaunay triangulation is most commonly used as it provides triangulations whose elements respect certain quality \bl{criteria}. \iffalse To describe the element quality for a 2D mesh, consider an equilateral triangle with area $A = \dfrac{A_{\Omega}}{N}$ as the standard element, where N is the number of elements in mesh and $A_{\Omega}$ refers to the area of the tessellated domain. \fi Given a set of points V, the Delaunay technique attempts to create triangulations wherein each triangle maximizes one (or more) of the following ratios: the inradius to the circumradius; the shortest edge to the longest edge; the shortest altitude to the longest edge; the aspect ratio, etc. \cite{parthasarathy1994comparison, field2000qualitative,  FIELD1991263,
dannelongue1991three, de1990explicit, parthasarathy1991constrained, opac-b1031284, BHATIA1990309, bank1996algorithm,watabayshi1986optimized, fukuda1972automatic, baker1989element, miller1997teng}. In this article, we use the inradius to circumradius ratio as our element-wise quality metric, given by: 
\begin{align}
Q = \frac{8A^2\bl{s}}{abc(a+b+c)}, 
\end{align}  
where a,b,c are side lengths, \bl{s} is the dimension and A is the area of the element. The use of the inradius to circumradius ratio for measuring the quality of elements was suggested by Cavendish, Field and Frey \cite{caendish1985apporach}. A high value of Q $\in [0,1]$ implies better quality \bl{elements}. An equilateral triangle and a standard tetrahedron has Q = 1. They are considered the standard elements for 2D and 3D meshes respectively. 
\section{Methods}
\label{sec:methods}

\subsection{Smoothing with the Shape Parameter}
\label{subsec:shapeP}

Our goal is to develop an iterative quasi-local smoothing algorithm to rectify singular deformation maps. To do so, we \bl{utilize an interesting feature} of RBF interpolation: smoothing using the shape parameter. This was first proposed by Beatson in the context of surface reconstruction from point cloud data~\cite{Carr:2003:SSR:604471.604495}, and has since been used as part of a numerical method for solving coupled PDEs~\cite{SWFKIJNMF2014}. This technique is very simple to apply: first, find the interpolation coefficients $\lambda_i(\epsilon^*)$, where $\epsilon^*$ is some small non-zero value. Then, when evaluating the interpolant, replace $\epsilon^*$ with $\epsilon$, where $\epsilon \neq \epsilon^*$. In other words, given an evaluation node set $X = \{\vx_j\}_{j=1}^M$, evaluate the interpolant at each point $\vx_j$ as
\begin{align}
\label{eq:rbfE1}
I_{\phi} f(\vx_j,\epsilon) = \sum_{i=1}^N \lambda_i(\epsilon^*) \phi(\epsilon,r_i(\vx_j)),
\end{align}
where $r_i(\vx_j) = \|\vx_j - \vx_i\|$. If $\epsilon < \epsilon^*$, this amounts to evaluating the coefficients against a slightly smoother basis than the one we interpolated with; this results in a smoothing; conversely, choosing $\epsilon > \epsilon^*$ can result in a sharpening of low-frequency details.

In this article, we present and utilize a simple generalization of the above approach: we allow $\epsilon$ to vary from point to point. In other words, we now evaluate the interpolant pointwise as
\begin{align}
\label{eq:rbfE2}
I_{\phi} f(\vx_j,\epsilon_j) = \sum_{i=1}^N \lambda_i(\epsilon^*) \phi(\epsilon_j,r_i(\vx_j)),
\end{align}
where $\epsilon_j >0, j=1,\hdots,M$ are now \emph{pointwise} shape parameters that potentially differ from the interpolation shape parameter $\epsilon^*$. Since $\phi$ has global support, this is still not entirely a local smoothing. However, compared to previous approaches which use a single $\epsilon$, our new approach constitutes a \emph{quasi-local} smoothing of the interpolant. In Section \ref{sec:algorithm}, we describe a technique which generates each $\epsilon_j$ given $\epsilon^*$, $X$ and $Y^d \subset Y^b$ samples of the deformation function on the boundary. \bl{Here, $Y^d = pY^b$, $0<p\le1$ and $Y^b$ is set of points on the boundary of deformed domain. The points on boundary are chosen based on the equation that describes the boundary and a boundary thickness parameter $\alpha$ set in step \ref{step:3} of the Algorithm \ref{Algorithm}. For example, if the 2D domain is a unit circle centered at origin, all the points that satisfy the equation $x^2 + y^2 = 1$ within a tolerance of $\alpha$ falls on the boundary. The parameter $p$ is chosen randomly and the subset is formed uniformly. The idea here is to show the efficacy of deformation map in deforming the entire domain even when we pick fewer points on the boundary.} As we will see in Section \ref{sec:experiments}, \bl{the scalar-valued RBF approximation and smoothing method described here, when applied in component-wise fashion to 2D and 3D problems,} gives intuitive results in the form of an appropriately smoothed set of output nodes $Y$.

\subsection{RBF-interpolation Based Iterative Algorithm for Mesh Generation and Quality Improvement}
\label{sec:algorithm}

In this section we present the RBF-interpolation based algorithm for generating curvilinear mesh and iterative smoothing, discuss implementation details of the same and provide a detailed analysis in terms of complexity.\\

\noindent \textbf{Overview} \\
\bl{Algorithm} \ref{Algorithm} describes the algorithm for obtaining a high-quality deformed mesh given a set of points in an initial undeformed domain and a set of parameters that control the deformation and smoothing process. Broadly, the procedure can be seen as a collection of following tasks:
\begin{enumerate}
\item Given the undeformed domain and samples of a deformation function on the boundary, interpolate the function (given by Equation \eqref{eq:rbfI}) to recover the deformation map in the interior of the domain.
\item \label{item:two} Tessellate the deformed domain and calculate element quality $Q = \underline{q}_e$ (see Section \ref{subsubsec:qualityC}).
\item Distribute the quality metric from the elements to vertices by averaging the quality of elements in 2-ring neighborhood around each vertex (see Section \ref{subsubsec:qualityC}).
\item For vertices with quality ($\underline{q}_v$) less than a predefined tolerance, reduce the shape parameter ($\epsilon$) by some factor (see Section \ref{subsubsec:modifyep}).
\item \label{item:penult} Evaluate the interpolant using the list of modified shape parameters ($\epsilon_j$) described by Equation \eqref{eq:rbfE2} to obtain an improved deformed mesh.
\item Repeat Steps \ref{item:two} through \ref{item:penult} until convergence defined by stopping criteria.
\end{enumerate}
\begin{algorithm}
\caption{RBF-based iterative algorithm for mesh generation and quality improvement}
\label{Algorithm}
\begin{algorithmic}[1]
\State Set $\delta \leftarrow$ scaling factor for shape parameter update term
\State Set $\sigma \leftarrow$  falloff of local Gaussian smoothing for shape parameter
\State \label{step:3}Set $\alpha \leftarrow$ thickness of boundary
\State $X_i \leftarrow$ $N_i \times s$ matrix of interior nodes on $\Omega(\mathbb{R}^s)$
\State $X_b \leftarrow$  $N_b \times s$ matrix of boundary nodes on $\partial \Omega(\mathbb{R}^s)$
\State $X = X_i \cup X_b $ $\leftarrow$ $N \times s$ matrix containing all nodes, $N = N_i + N_b \leftarrow |\Omega|$
\State $X^d \subset X_b$, $N_d \times s$ matrix of data sites on $\partial \Omega$, $N_d \subset N_b = |\partial \Omega|$
\State $Y^d$ $\leftarrow$ $N_d \times s$ matrix of deformed boundary nodes (corresponding to $X^d$)
%\State Set $\rho \leftarrow$ tolerance for desired quality to flag nodes for correction
\State  $\kappa_t$ $\leftarrow$ desired target condition number of interpolated matrix 
\State \label{epIdeal} $\epsilon^*$ $\leftarrow$ ideal shape parameter corresponding to $\kappa_{\bl{t}}$
\State Initialize $\underline{\epsilon} = \epsilon^* \forall$ $\vx$ in $X$ 
\State $A$ $\leftarrow$ RBF interpolation matrix using $\epsilon^*$ 
\State \underline{$\lambda$} $\leftarrow$ Interpolation coefficients, obtained formally by finding $A^{-1}X^d$ once
\State  $Y$  $\leftarrow$ Evaluate the RBF interpolant built on \bl{$X^d$} at all nodes in $X$ with $\underline{\epsilon}$
\State Tessellate $Y$ to obtain element set $E$
\State For each $\vy$ in $Y$, store its $n_k$ 2-ring neighbors $\leftarrow \{\vy,\vy_k\}$
\For{each correction iteration until convergence}
\State \label{qCal} Calculate $\underline{q}_e \leftarrow$ quality per element in $E$
\State \label{hUpdate} \bl{Append $\|q_e\|_2$ to $\underline{h}$, history of mesh quality over iterations}
\State \label{cCheck} Check convergence: $\|\underline{q}_e\|_2<$  max(h)
\State \label{qPerV} Distribute $\underline{q}_e$ to $\underline{q}_y \leftarrow$ quality per vertex in $Y$
\For{each point $\vy_k$ in $Y$}
%\State $\underline{\Psi}$ $\leftarrow q_y(\{\vy^p,\vy^p_k\})$, quality of $\vy$ and $n_k$ 2-ring neighbors 
%\If { min($\underline{\Psi}$) $< \rho $ }  
\State $\mu(\vy_k,\alpha) \leftarrow \min \|\vy_k-\vy_{p}\|_{2}$ $\forall$ $\vy_p \in Y^d$ 
\State $\mu(\vy_k,\alpha) = 0$ if $\mu(\vy_k, \alpha) \leq \alpha$ 
\For{ each $j$ from 1 to $n_k$ }
\State $\lf(\underline{\Psi}_k\rt)_j \leftarrow |q_k - q_j|$
\EndFor
\State $\underline{\gamma}_k \leftarrow e^{-\sigma \underline{\Psi}_k}$
\State $\theta_k \leftarrow \delta \mu(\vy_k,\alpha)$
\State \label{epUpdate}   $\underline{\epsilon}_k \leftarrow  \underline{\epsilon}_k - \theta_k \underline{\gamma}_k$
%\EndIf
\EndFor
\State \label{yCal} Compute smoothed node set $Y$ using new $\underline{\epsilon}$ and precalculated $\underline{\lambda}$
\State \label{eCal} Tessellate $Y$ to obtain element set $E$
\State For each $\vy$ in $Y$, update its $n_k$ 2-ring neighbor stencil $\leftarrow \{\vy,\vy_k\}$
\EndFor
\end{algorithmic}
\end{algorithm}
\subsubsection{Computing quality per-element and per-vertex} 
\label{subsubsec:qualityC}
At each iteration of our RBF-based technique, the resulting mesh element quality ($\underline{q}_e$) is determined based on one of the definition of quality metric as detailed in Section \ref{subsec:MeshQ}. The overall quality of the mesh is the aggregate of quality of all elements in \bl{the} mesh. This is used to determine the stopping criteria for the algorithm. If the overall quality satisfies a predefined threshold, the algorithm converges. 

Let $\vy_k$ be the 2-ring neighbors of a node $\vy$ in $Y$. We view these nodes as constituents of a stencil for measuring \bl{the} per-vertex quality. Here, the number of vertices in the stencil ($n_k$) depends on the degree of connectedness of the vertex. For instance, vertices which are close to the domain boundary have fewer neighbors while others have a full connectivity with 2-ring neighbors. \bl{Because the quality is defined per-element and we want to have a quality measure per-vertex, we need to find the elements connected by a vertex}. To aggregate elemental qualities $\underline{q}_e$ to individual vertex qualities $\underline{q}_y$ for \bl{vertices} in $Y$, we use the following average:
\begin{align}
\label{eq:qEtoqV}
q(\vy_k) = (q_y)_k = \frac{1}{n_k}\sum_{i = 1}^{n_k} {q}_e (\vy_i), k = 1,\hdots,|Y|,
\end{align}
where $|Y|$ is the total number of vertices in the domain.

\subsubsection{Modifying shape parameter based on \bl{the} per-vertex quality}
\label{subsubsec:modifyep}

\bl{We now describe our formula for generating a new modified shape parameter at each vertex. We modify the shape parameter at a vertex based on two factors: the quality measure at the vertex (given by Equation \eqref{eq:qEtoqV}), and the proximity of the vertex to the boundary.  Without loss of generality, we focus on the vertex $\vy_k$. Let $\underline{\epsilon}_k^{old}$ be the $n_k$-long vector of shape parameters of $\vy_k$ and its 2-ring neighbors in the current iteration. At the first iteration of the smoothing algorithm, the shape parameters at all vertices are the same, \emph{i.e.}, $\underline{\epsilon}^{old}_k = \epsilon^*, k=1,\hdots,N$. The goal is to obtain $\underline{\epsilon}_k^{new}$, the new vector of shape parameters, for every subsequent iteration. We propose a simple update of the form
\begin{align}
\label{eq:upEp}
\underline{\epsilon}^{new}_k = \underline{\epsilon}^{old}_k - \theta_k \underline{\gamma}_k,
\end{align}
where $\theta_k$ is a factor that accounts for proximity to boundaries, and $\underline{\gamma}_k$ is a factor that depends on the vertex qualities of $\vy_k$ \emph{and} its 2-ring neighbors; this formula is given on line 30 of Algorithm 1. 
We will first explain the $\underline{\gamma}_k$ term, then the $\theta_k$ term.}

\bl{The term $\underline{\gamma}_k$ is a function of the vertex quality $(q_y)_k$ associated with the vertex $\vy_k$. Specifically, this term attempts to decrease $\underline{\epsilon}^{old}_k$ whenever the vertex quality associated with $\vy_k$ is significantly different from the vertex qualities of its 2-ring neighbors. First, we define the quantity $\underline{\Psi}_k$ as
\begin{align}
\lf(\underline{\Psi}_k\rt)_j = |(q_y)_k - q_j|, j=1,\hdots,n_k,
\end{align}
where $j$ indexes the 2-ring neighbors of the vertex $\vy_k$. Clearly,$\underline{\Psi}_k$ is a vector of differences in quality between $\vy_k$ and its 2-ring neighbors. Our formula for $\underline{\gamma}_k$ satisfies two requirements: first, that $\underline{\gamma}_k$ change \emph{smoothly} as a function of $\underline{\Psi}_k$, and second, that $\underline{\gamma}_k$ is smaller as we go further away from $\vy_k$. These two requirements are satisfied by requiring $\underline{\gamma}_k$ to take the form
\begin{align}
\underline{\gamma}_k = e^{-\sigma \underline{\Psi}_k},
\end{align}
where $\sigma$ is some user-supplied \emph{falloff} factor. If $\sigma$ is small relative to the distance between nodes, the different $\underline{\Psi}_k$ values contribute more equally to $\underline{\gamma}_k$. In contrast, if $\sigma$ is large, the contributions of $\underline{\Psi}_k$ corresponding to nodes other than $\vy_k$ are smaller. In this article, we use values of $\sigma$ that ensure that we are in the latter regime. This allows us to more effectively correct localized irregularities in vertex quality, while still smoothly updating the $\underline{\epsilon}_k^{old}$ values.}

\bl{When we attempted to update $\underline{\epsilon}_k^{old}$ using only the $\underline{\gamma}_k$ values, we ran into two difficulties. First, we observed that nodes from the interior would leave the domain boundary, and hence would need to be periodically deleted from the domain. Second, such updates tended to undo mesh refinement near the domain boundary. Our first attempt at fixing this problem was to multiply $\underline{\gamma}_k$ by a \emph{switch} that \emph{turns off} smoothing near the boundary. However, noticing that this produced some mesh tangling outside the boundary-refined layers, we choose instead to multiply $\underline{\gamma}_k$ by the scalar term $\theta_k$ (line 29 in Algorithm 1) defined as
\begin{align}
\theta_k = \delta \mu(\vy_k,\alpha).
\end{align}
Here, $\delta$ is some small number that controls the magnitude of $\theta_k$, and $\mu(\vy_k,\alpha)$ is a function that effectively specifies a ``boundary-layer'' for our algorithm. Let $\vy_p$ be the closest boundary point to $\vy_k$. Then, $\mu(\vy_k,\alpha)$ (lines 23 and 24 of Algorithm 1) is defined as
\[
\mu(\vy_k,\alpha) = 
\begin{cases}
\|\vy_k - \vy_p\|, & \|\vy_k - \vy_p\| > \alpha\\
0, & \|\vy_k - \vy_p\| \leq \alpha
\end{cases}
\]
This function ensures that no update is made to the shape parameter of any node $\vy_k$ within distance of $\alpha$ from its closest point $\vy_p$ on the boundary. Further, nodes $\vy_k$ further away from their closest boundary points are allowed to receive larger updates to their shape parameter vectors $\underline{\epsilon}_k$.}

\bl{In general, we find that $\delta$ needs to be small to improve quality, ensuring that the shape parameters are not decreased too much in any iteration. Currently, $\delta$, $\alpha$ and $\sigma$ are selected by trial and error, but one could imagine using training techniques from the neural networks literature to accomplish this. We leave such extensions for future work.}

\subsubsection{Stopping criterion}
We now present a stopping criterion for the smoothing algorithm. The criterion is designed to stop the iterative smoothing if the mesh quality begins to worsen as a consequence of the iterative procedure. \bl{Such a worsening in quality, when it occurs, is a consequence of the global support of the RBF interpolant. Despite the local nature of the shape parameter updates, the global support of the RBFs means that most nodes are moved to some extent.} 

To determine a good stopping point for the iterative smoothing process, we simply check the 2-norm of the per-vertex quality measure, \emph{i.e.}, $\|\underline{q}_e \|_2$. If $\|\underline{q}_e \|_2$ is \bl{smaller} for the current iteration than for previous ones, the algorithm halts. This choice of stopping criterion may not be ideal, since it aggregates mesh quality into a single number. However, we have found that it works well in conjunction with the global RBF interpolant. We leave the question of stopping criteria for future work.

\section{Complexity Analysis}
\label{sec:complex}

\subsection{Preprocessing}
We first consider the preprocessing costs of our algorithm. Consider a tessellated domain $\Omega \subset \mathbb{R}^s$. Let $N_b$ be the number of points on the boundary of \bl{the} domain ($\partial \Omega$) and $N_i$ be the points in the interior. The total number of points in the domain is then given by $N = N_i + N_b$. However, not all $N_b$ points are used to recover the deformation map via the RBF interpolant. Let $N_d \subset N_b$ be the number of points used to build the RBF interpolant. The initial preprocessing step involves computing and decomposing the RBF interpolation matrix once for a cost of $O(N_d^3)$. The interpolant can then be evaluated for $O(N N_d)$. However, it is more intuitive to express this cost in terms of the number of interior points. We now present that derivation, specialized to $s=2,3$.

\subsubsection{Complexity analysis in 2D ($s=2$)}
\bl{
Before proceeding, assume that points on the boundary are evenly-spaced with spacing $h_b$. \bl{Further, assume that the interior nodes are spacing $h_i$. Then, we have}
\begin{align}
\label{eq:hout}
h_b = \frac{l}{N_b},  h_i = \lf(\frac{A}{N_i}\rt)^{\frac{1}{2}},
\end{align} 
where $l$ is the perimeter of the boundary and \bl{$A$ is the area of the domain}. Assuming without loss of generality that $h_i = h_b$, we have  
\begin{align}
\frac{l}{N_b} &= \lf(\frac{A}{N_i}\rt)^{\frac{1}{2}} 
\implies N_b = l A^{-\frac{1}{2}} N_i^{\frac{1}{2}}. \label{eq:relation_nb_ni_A}
\end{align}
Now, letting $N_d = p N_b, 0<p\leq 1$, we can rewrite Equation \eqref{eq:relation_nb_ni_A} as
\begin{align}
\label{eq:relation_nd_nb}
N_d = p N_b = pl A^{-\frac{1}{2}} N_i^{\frac{1}{2}}.
\end{align}
Since the interpolation matrix is inverted for a one-time cost of $O(N_d^3)$, we now have an explicit expression for that cost. Using Equation \eqref{eq:relation_nd_nb}, this cost becomes:
\begin{align}
N_d^3 = p^3 l^3 A^{-\frac{3}{2}} N_i^{\frac{3}{2}} 
\implies N_d^3 = p^3 A^{\frac{1}{2}} \lf(\frac{l}{A} \rt)^3 N_i^{\frac{3}{2}}. \label{eq:gc_pre}
\end{align}
Let $\psi_s$ be a domain-dependent constant in $s$ dimensions so that for $s=2$, $\psi_2 = \frac{l}{A}$. We can use this to rewrite
Equation \eqref{eq:gc_pre}, obtaining: 
\begin{align}
\label{eq:penult_eq}
N_d^3 = p^3  A^{\frac{1}{2}} \psi_2^3 {N_i^{\frac{3}{2}}}.
\end{align}
In general, $N_b << N_i$, implying that $N \approx N_i$. Thus, the preprocessing cost $C_2$ for our technique in $s=2$ spatial dimensions is asymptotically $C_2 = O(N^{1.5})$. Note that if the interpolation problem could be solved in $O(N_d)$ operations, this cost reduces to $O(\sqrt{N})$. We leave this extension for future work.}

\subsubsection{Complexity analysis in 3D ($s=3$)}
\bl{
We now derive 3D complexity estimates for our preprocessing step. Assuming that nodes on the domain boundary (now a surface of co-dimension one in $\mathbb{R}^3$) are quasi-uniform with spacing $h_b$ \bl{and assuming the interior node spacing is $h_i$, we have}
\begin{align}
\label{eq:hin_3D}
h_b = \lf( \frac{A}{N_b} \rt)^{\frac{1}{2}} \textrm{ and } h_i = \lf(\frac{{V}}{{N_i}}\rt)^{\frac{1}{3}},
\end{align} 
where $A$ is now the surface area of $\partial \Omega$, and $V$ is the volume of $\Omega$. Assuming again that $h_i = h_b$, we have
\begin{align}
\lf(\frac{A}{N_b}\rt)^{\frac{1}{2}} &= \lf(\frac{{V}}{{N_i}}\rt)^{\frac{1}{3}} \implies N_b = AV^{-\frac{2}{3}} N_i^{\frac{2}{3}}. \label{eq:relation_nb_ni_A_V}
\end{align}
Expressing this in terms of $N_d$, the number of data sites used \bl{to} build the interpolant, we have
\begin{align}
\label{eq:relation_nd_nb_3D}
N_d = p N_b = p AV^{-\frac{2}{3}} N_i^{\frac{2}{3}}.
\end{align}
The preprocessing cost is $O(N_d^3)$, which is now given by:
\begin{align}
N_d^3 =  p^3 A^3 V^{-2} N_i^2.
\end{align}
Now letting $\psi_3 = \frac{A}{V}$, we have
\begin{align}
N_d^3 =  p^3 V \psi_3^3 N_i^2.
\end{align}
The preprocessing cost $C_3$ can now be expressed in terms of the total number of points $N$ as $C_3 = O(N^{2})$. Again, as in 2D, it is possible to lower this cost (to $O(N^{\frac{2}{3}})$) if the interpolation problem is solved in $O(N_d)$ operations.}

\subsubsection{Finding the initial shape parameter}
Another contribution to the preprocessing cost comes from the calculation of the initial shape parameter ($\epsilon^*$). We use the fzero function in Matlab to find this shape parameter. This function uses an iterative method called the Brent-Dekker method to find the zero of a function in a given interval. Consequently, it requires an evaluation of that function multiple times. In our application, the function that must be evaluated is the condition number of the RBF interpolation matrix. This can be computed for a cost of $O(N_d^3)$ if the 2-norm condition number is used, and a cost of $O(N_d^2)$ if the 1-norm or max-norm condition numbers are used. From Equation \eqref{eq:relation_nd_nb}, \bl{and considering max-norm}, it is obvious that this cost scales as $O(N)$ in 2D; similarly, Equation \eqref{eq:relation_nd_nb_3D} \bl{for max-norm condition number} shows that this cost scales as $O(N^{\frac{4}{3}})$ in 3D.

\subsection{Complexity of the Smoothing Algorithm}
We now analyze the complexity of a single step of our smoothing algorithm. To do so, we break down our algorithm into several key steps.
 
\subsubsection{Finding the 2-ring neighbors}
Each iteration of the algorithm requires finding the 2-ring neighbors of each vertex using the Delaunay triangulation (which is itself constantly being updated). To find the 2-ring neighbors of each vertex, we first find the list of vertices connected to each vertex (the 1-ring neighbors). The cost of this operation for N vertices scales as $O(N \log N)$, with a dimension dependent constant. The next step is to find the set of immediate neighbors of the 1-ring neighbors. To do so, we simply repeat the above step for each of the 1-ring neighbors. The total asymptotic cost of finding the 2-ring neighbors is therefore:
\begin{align}
C_{\textrm{$n_k$-neighbors}} = N(1 + n_k)\log N.
\end{align}

\subsubsection{Calculating per-element and per-vertex quality}
At each iteration, we calculate per-element quality ($\underline{q}_e$) and distribute it to the constituent nodes forming the elements as ($\underline{q}_y$) given by Equation \eqref{eq:qEtoqV}. $\underline{q}_e$ can be computed for a cost of $O(N)$. Similarly, $\underline{q}_y$ requires an averaging over elements connected to each vertex. If the average number of elements connected to each vertex is $n$, then this cost scales as $O(nN)$, where $n << N$. The complexity of this step therefore scales as $O(N)$.

\subsubsection{Updating the evaluation shape parameter}
At each iteration, the algorithm updates the shape parameter for each vertex in the \bl{node set} based on the quality metric and calculations described by
Equation \eqref{eq:upEp}. This operation utilizes the 2-ring neighbor information from previous step and the predefined parameters described in Section \ref{subsubsec:modifyep}. By a similar argument to the previous subsection, this update also scales as $O(nN) \approx O(N)$.

\subsubsection{Computing the smoothed \bl{node set} }
To obtain the smoothed \bl{node set} $Y$ at each iteration, we need to compute the RBF evaluation matrix and multiply it with the precomputed interpolation coefficients. The computation of the evaluation matrix is straightforward, as shown in line \ref{yCal} of \bl{Algorithm} \ref{Algorithm}. The operation scales as $O(N N_d)$ when $N$ evaluation points are used. It is clear from Equation \eqref{eq:relation_nd_nb} that this cost scales as $O(N^{\frac{3}{2}})$ in 2D with a small constant term. In 3D, Equation \eqref{eq:relation_nd_nb_3D} shows that this cost scales as $O(N^{\frac{5}{3}})$, again with a small constant. 
%MIKE
%In practice, we have found that these costs are subdominant compared to the cost of per-iteration tessellation.

\subsubsection{Tessellation of domain to obtain element set} 
As a last step of each smoothing iteration, the \bl{node set} $Y$ is tessellated to obtain a mesh which is smoother than the previous iteration. This operation is performed using the Delaunay triangulation in the code (refer to line \ref{eCal} in \bl{Algorithm} \ref{Algorithm}). There are many other  algorithms to obtain a mesh from the \bl{node set} and depending on the use case, this choice can vary. In general, this cost is $O(N \log N)$ in 2D, and $O(N^2)$ in 3D.
\section{Results} \label{sec:experiments}

We now present the results of our numerical experiments using our algorithm. To demonstrate the efficacy of the technique across different types of problems, we focus on three test cases involving domains with boundaries of different smoothness:
\begin{enumerate}
\item Deforming a $C^1$ boundary to a $C^\infty$ boundary.
\item Deforming a domain with a $C^\infty$ outer boundary and a $C^\infty$ inner boundary (an annulus) to a $C^1$ outer boundary and $C^1$ inner boundary (a square with an airfoil cavity).
\item Deforming a cube ($C^1$ boundary) to a sphere ($C^{\infty}$ boundary).
\end{enumerate}
A fourth category for which results are not shown is the deformation of domains with $C^{\infty}$ boundaries to domains with $C^{\infty}$ \bl{boundaries}. We do not show results for this case, as the smoothing procedure is completely unnecessary here. The interpolation step itself produces excellent meshes, at least partly due to the spectral convergence rates achieved by the RBF interpolant to the deformation map.

In order to obtain tessellations on the undeformed domains, \bl{we first generate a set of reasonably well-distributed nodes using a repulsion algorithm such as the one used by Distmesh\cite{perssondistmesh}. This gives us both interior and boundary nodes. The undeformed mesh is then generated by applying a simple Delaunay triangulation on this set of nodes.} Following Section \ref{subsec:shapeP}, we computed initial shape parameters for all the above test cases. The results are summarized in Table \ref{tab:TestParams}. 
\begin{table}[H]
\caption{Summary of test parameters. $C^4$ Mat\'ern kernel was used for all the experiments}
\label{tab:TestParams}

\begin{tabular}{>{\raggedright\arraybackslash}m{0.2\textwidth} >{\centering\arraybackslash}m{0.03\textwidth} >{\centering\arraybackslash}m{0.17\textwidth} >{\centering\arraybackslash}m{0.04\textwidth} >{\centering\arraybackslash}m{0.07\textwidth} >{\centering\arraybackslash}m{0.04\textwidth} >{\centering\arraybackslash}m{0.12\textwidth} >{\raggedright\arraybackslash}m{0.05\textwidth}}

\hline\noalign{\smallskip}
\multicolumn{1}{>{\centering\arraybackslash}m{0.2\textwidth}}{Test} & \multicolumn{1}{>{\centering\arraybackslash}m{0.03\textwidth}}{ $\epsilon^*$ } & \multicolumn{1}{>{\centering\arraybackslash}m{0.17\textwidth}}{$N = N_i + N_b$}& \multicolumn{1}{>{\centering\arraybackslash}m{0.04\textwidth}}{$N_d$}& \multicolumn{1}{>{\centering\arraybackslash}m{0.07\textwidth}}{Mesh norm }&
\multicolumn{1}{>{\centering\arraybackslash}m{0.04\textwidth}}{$\alpha$}&
\multicolumn{1}{>{\centering\arraybackslash}m{0.12\textwidth}}{$\delta$}&
\multicolumn{1}{>{\centering\arraybackslash}m{0.05\textwidth}}{$\sigma$}\\

\noalign{\smallskip}\hline\noalign{\smallskip}
$C^1$ boundary to a $C^\infty$ boundary & 0.2497 & 2106 $=$ 1906$+$200 & 172 & 1.0058 & 0.001 & 9.9422e-07 & 0.1006\\
Annulus to a square with an airfoil cavity & 0.4556 &  1420 $=$ 1220$+$200 & 146 & 0.0137 & 0.01 & 7.3014e-06 & 1.3696\\
Cube to a sphere & 0.7354 & 20063 $=$ 18563$+$1500 & 1304 & 0.0015 & 0.001 & 6.5200e-07 & 0.0153\\
\noalign{\smallskip}\hline
\end{tabular}

\end{table}
For all experiments, our general approach is to prescribe a boundary deformation at a \emph{subset} of the boundary nodes using a conformal map (or some other transformation). The RBF \bl{interpolant} to the deformed boundary is then used as a proxy for the deformation map, and evaluation in the interior of the deformed domain gives us a set of interior points on that domain. We then run a few steps of our iterative smoothing algorithm to improve our mesh. \bl{In order to compare the RBF based smoothing with other smoothing techniques, we consider the Laplace smoother. At the end of each experiment in 2D, we run the same number of Laplace smoothing iterations as required by RBF-based smoothing before the stopping criterion terminates the algorithm. We also document the computational cost for each step of Algorithm 1 in Table \ref{tab:runtime} for all our experiments.}
\begin{table}[h!]
% table caption is above the table
\caption{A summary of runtime (in seconds) for all the experiments. The table shows time-averaged per-iteration costs for the important steps in Algorithm 1 (from Section \ref{sec:algorithm}) for each of our experiments.}
\label{tab:runtime}      % Give a unique label
% For LaTeX tables use

\begin{tabular}{>{\raggedright\arraybackslash}m{0.1\textwidth} >{\centering\arraybackslash}m{0.25\textwidth} >{\centering\arraybackslash}m{0.25\textwidth} >{\centering\arraybackslash}m{0.25\textwidth}}
\hline\noalign{\smallskip}
\multicolumn{1}{>{\centering\arraybackslash}m{0.1\textwidth}}{Step No.} & \multicolumn{1}{>{\centering\arraybackslash}m{0.25\textwidth}}{ $C^1$ boundary to $C^\infty$ boundary} & \multicolumn{1}{>{\centering\arraybackslash}m{0.25\textwidth}}{Annulus to a square with an airfoil cavity}& \multicolumn{1}{>{\centering\arraybackslash}m{0.25\textwidth}}{Cube to a sphere}\\
\noalign{\smallskip}\hline\noalign{\smallskip}
12 & 0.000621 & 0.000612 & 0.301047 \\ 
 %\hline
 13 & 0.000613 & 0.000334 & 0.298201\\
 %\hline
 14 & 0.002114 & 0.002389 & 0.746076 \\
 %\hline
 15 & 0.004944 & 0.009064 & 0.226799 \\
 %\hline
 16 & 2.722284 & 1.563437 & 386.769592\\
 %\hline
 18 & 0.000168 & 0.000093 & 14.133027 \\
 %\hline
 19 & 0.000189 & 0.000011 & 0.000015\\
 %\hline
 20 & 0.000238 & 0.000018 & 0.000011\\
 %\hline
 21 & 0.235500 & 0.382419 & 0.541982\\
 %\hline
 23-30 & 0.303988 &  0.118030 & 13.342620 \\
 %\hline
 32 & 0.059021 & 0.152558 & 5.322097\\
 %\hline
 33 & 0.005053 & 0.020471 & 0.005053\\
 %\hline
 34 & 1.751741 & 3.605688 & 384.959461\\
 %\hline
\noalign{\smallskip}\hline
\end{tabular}

\end{table}
\subsection{Deforming a $C^1$ boundary to a $C^\infty$ boundary}
\label{sec:cat1}
\begin{figure}[h!]
\centering
\subfigure[Undeformed domain]{%
 \includegraphics[width=.43\textwidth]{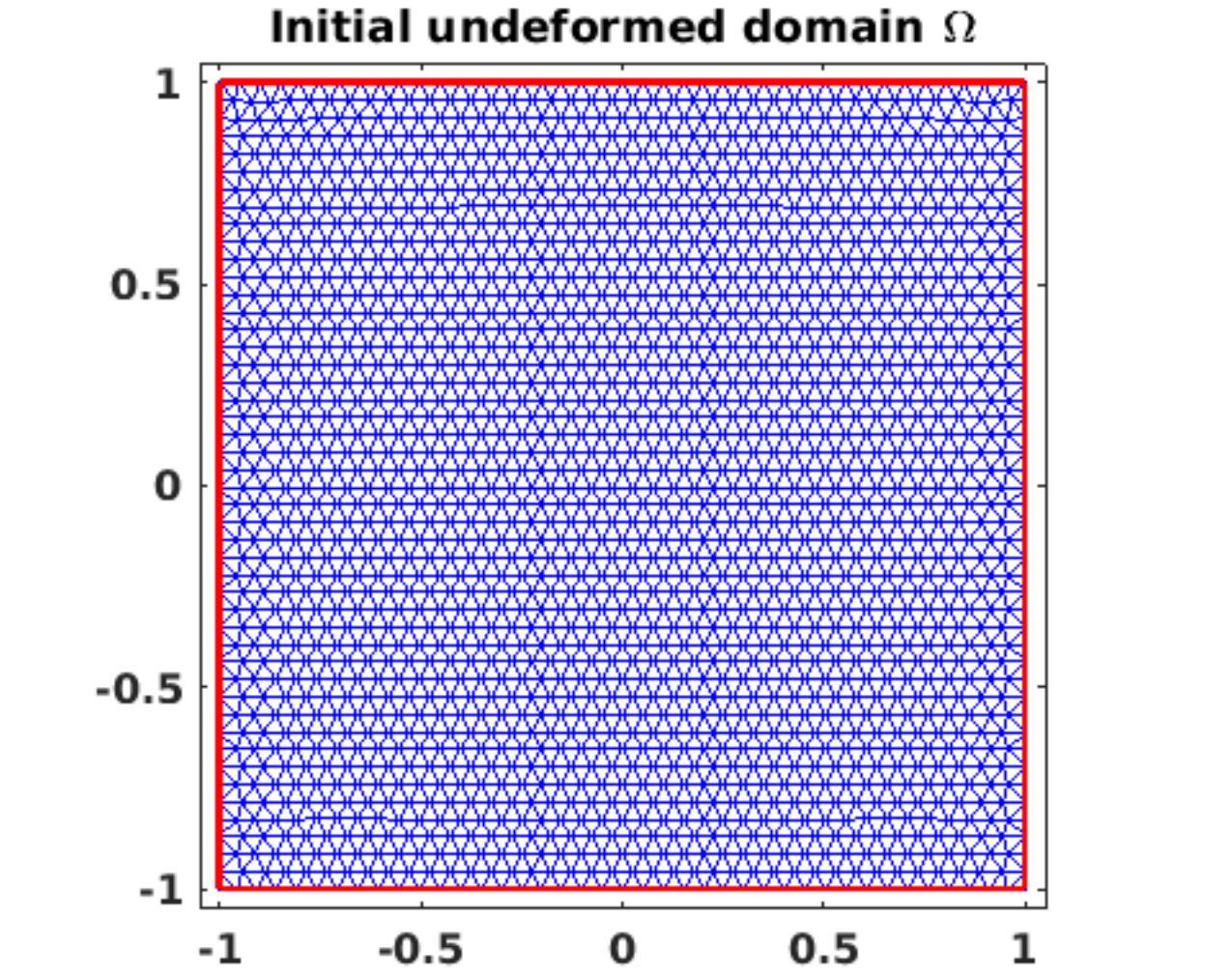}
\label{fig:1a}}
\subfigure[Quality of undeformed mesh]{%
 \includegraphics[width=.43\textwidth]{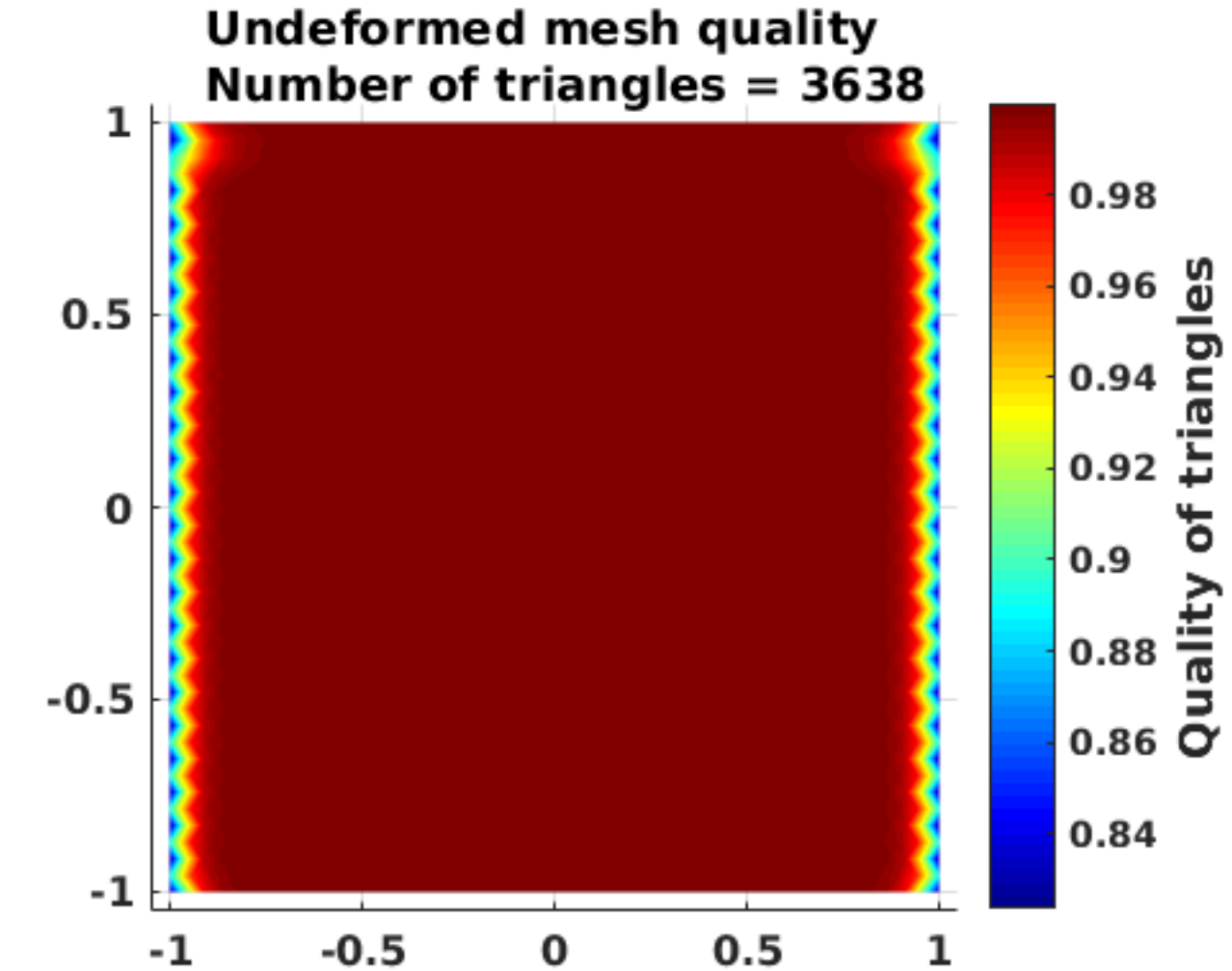}
\label{fig:1b}}

\subfigure[Deformed mesh before smoothing]{%
 \includegraphics[width=.43\textwidth]{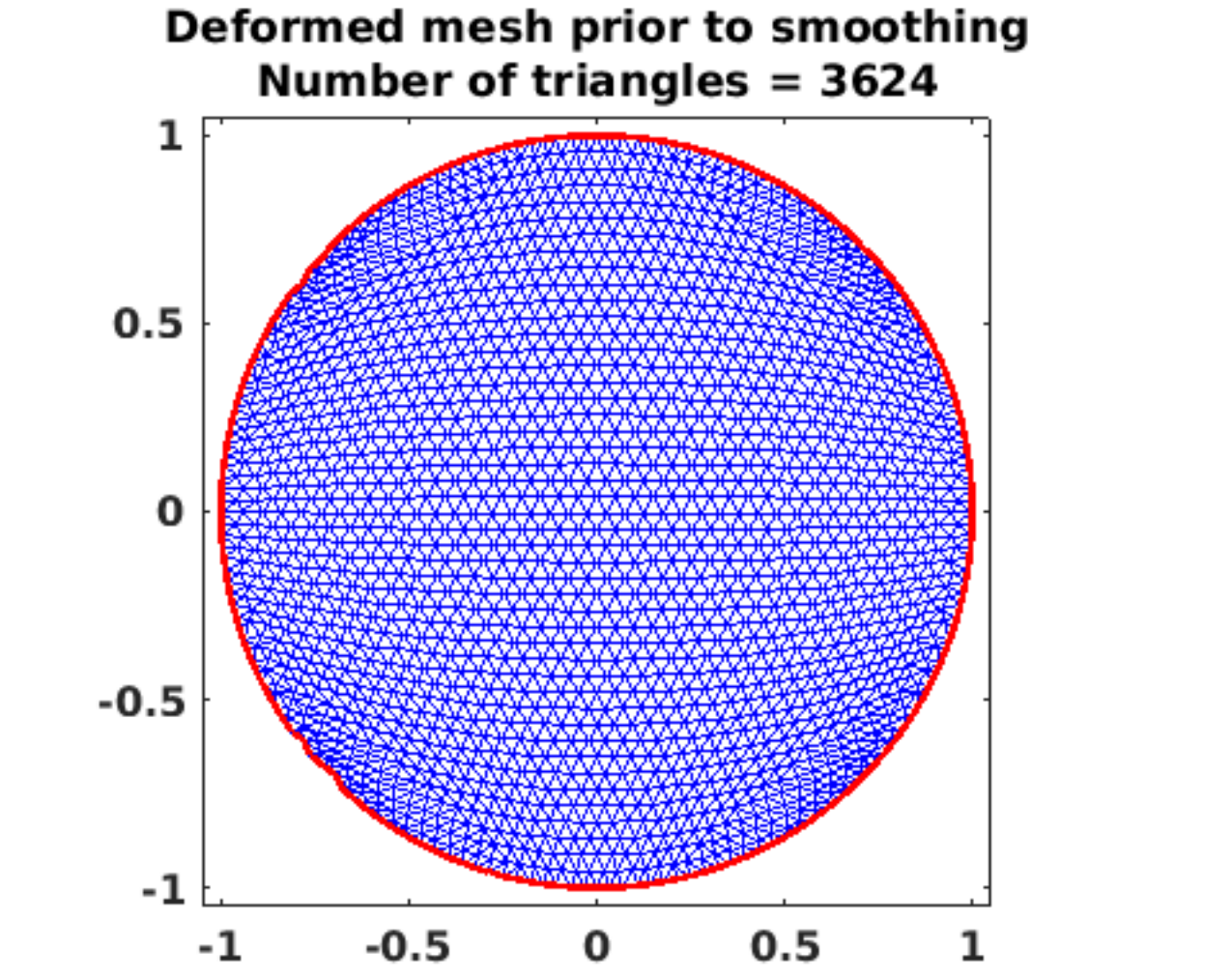}
\label{fig:1c}}
\subfigure[Quality of deformed mesh before smoothing]{%
 \includegraphics[width=.43\textwidth]{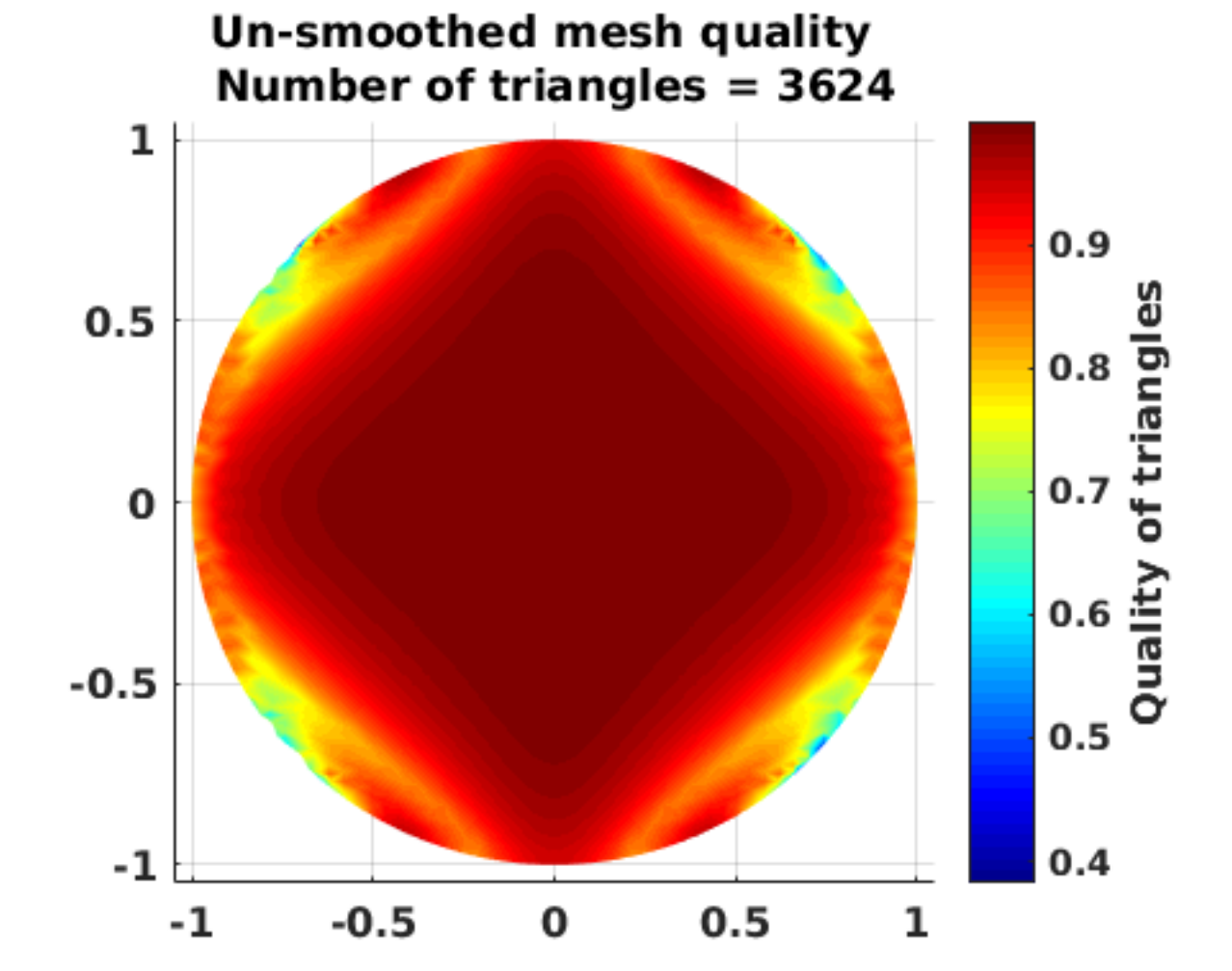}
\label{fig:1d}}

\subfigure[Deformed mesh after smoothing]{%
 \includegraphics[width=.43\textwidth]{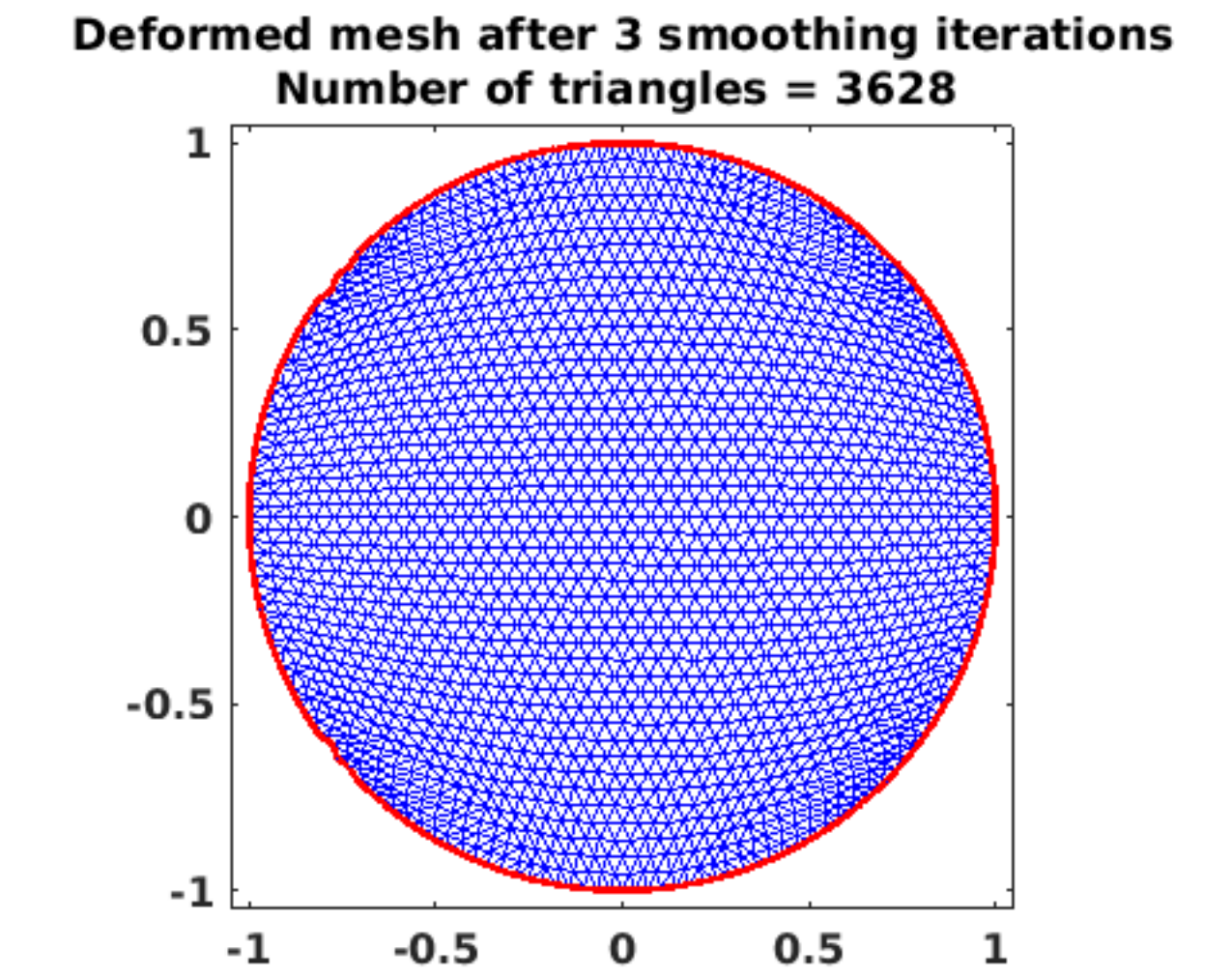}
\label{fig:1e}}
\subfigure[Quality of deformed mesh after smoothing]{%
 \includegraphics[width=.43\textwidth]{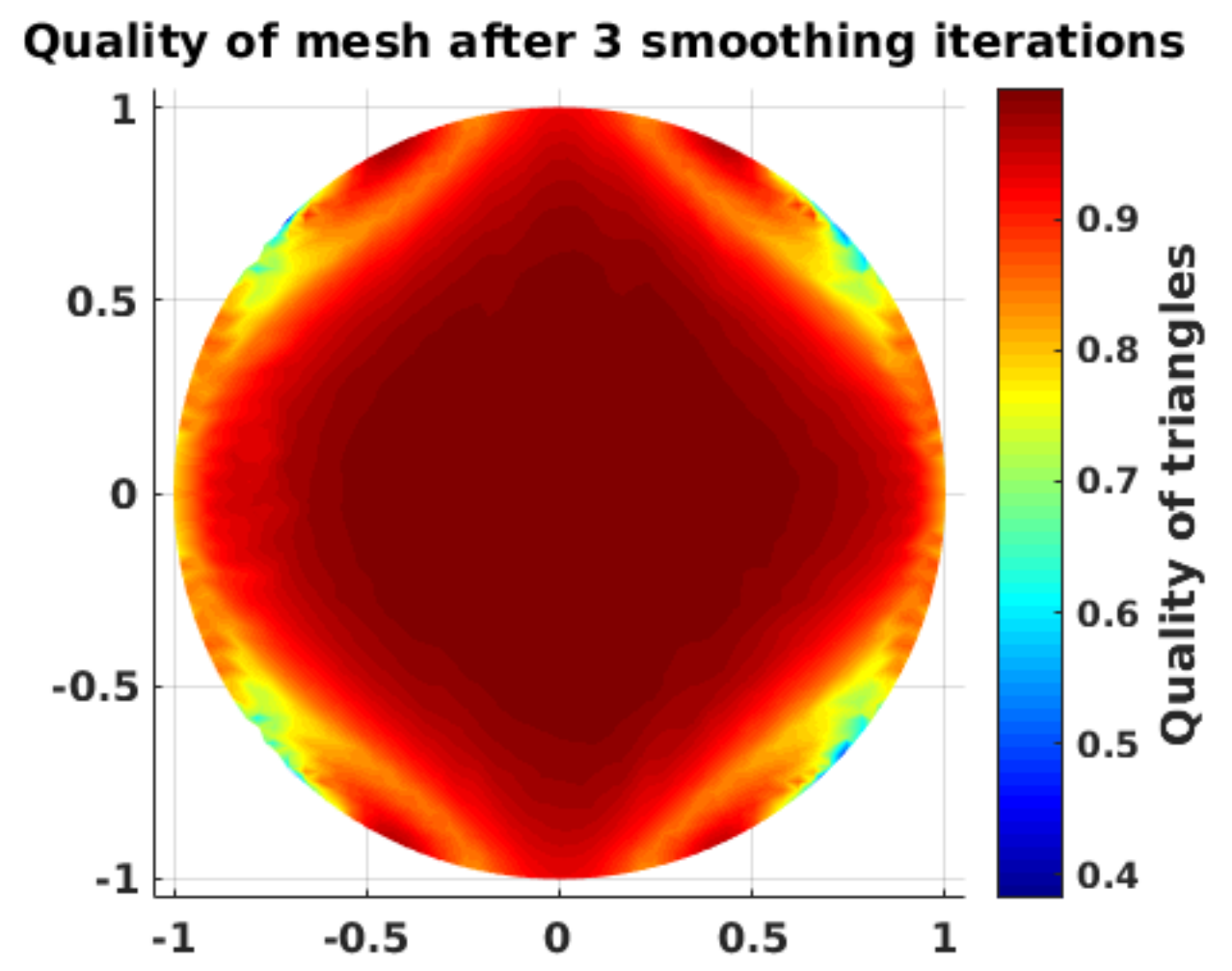}
\label{fig:1f}}
\caption{Deforming a $C^{1}$ boundary to a $C^{\infty}$ boundary.}
\label{fig:figuretype1a}  
\end{figure} 
For this test, we set $\Omega$ to be the unit square $[0,1]^2$; naturally, its boundary is of limited smoothness, \emph{i.e.}, $\partial \Omega \in C^1(\mathbb{R}^2)$. We then prescribe a deformation so that the deformed boundary $\partial \Omega'$ is a circle, and $\partial\Omega' \in C^\infty(\mathbb{R}^2)$. The deformation map $f$ is given component-wise by:
\begin{align}
f(x,y) = \lf[x \lf({1-\frac{y^2}{2}} \rt)^{\frac{1}{2}}, y \lf({1-\frac{x^2}{2}} \rt)^{\frac{1}{2}}\rt].
\end{align}
This conformal map has four singularities corresponding to the four corners of the square. These singularities are likely to manifest as distortions in the tessellation of the deformed domain. The results of this experiment are shown in Figure \ref{fig:figuretype1a}.
\begin{figure}[h!]
\centering
\includegraphics[width=.43\textwidth]{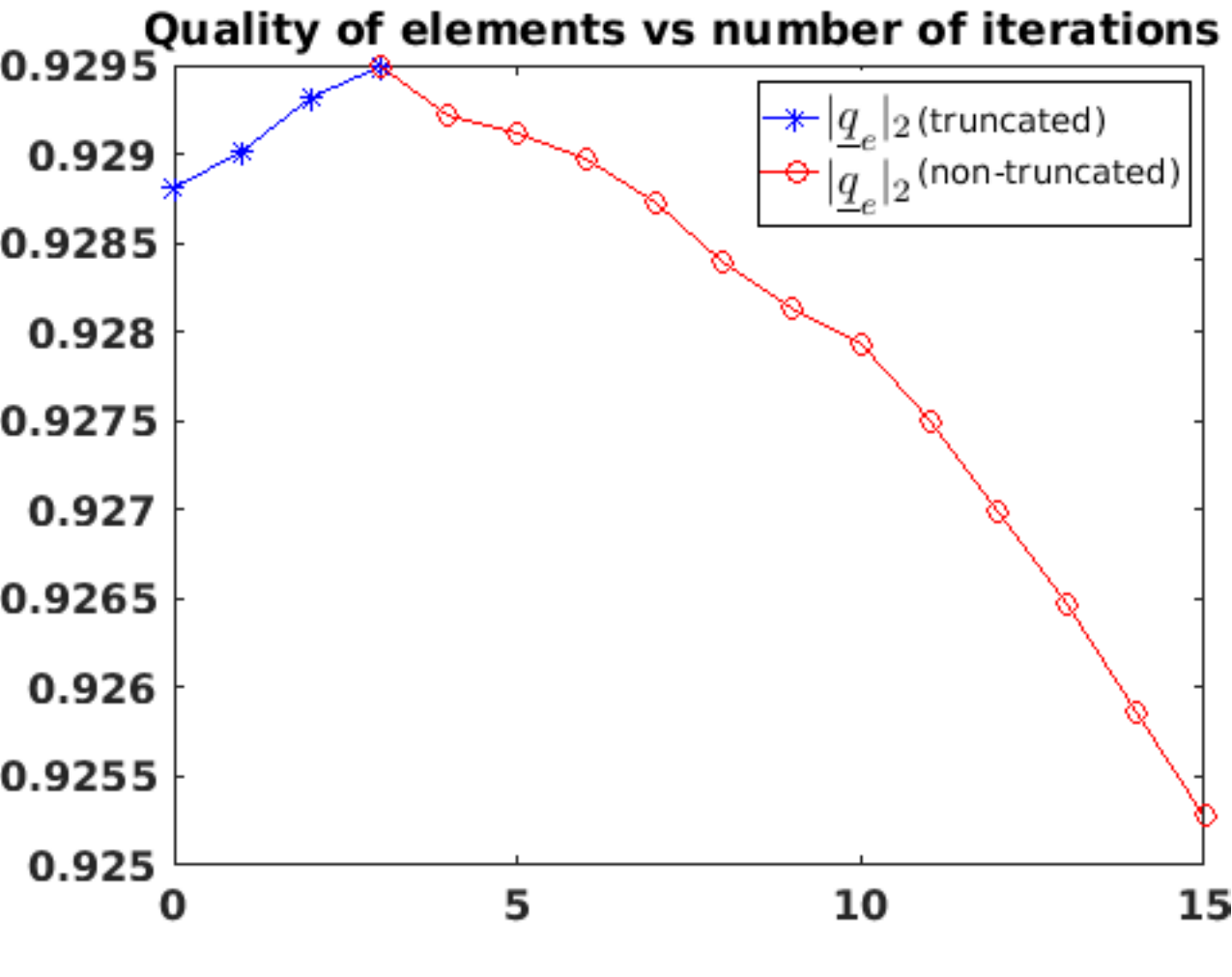}
\caption{Quality as a function of number of iterations when deforming a $C^{1}$ boundary to a $C^{\infty}$ boundary. The blue line corresponds to the use of a stopping criterion. The red line shows the quality when no stopping criterion is used.}
\label{fig:figuretype1c}  
\end{figure} 
\begin{figure}[h!]
\centering
\includegraphics[width=.50\textwidth]{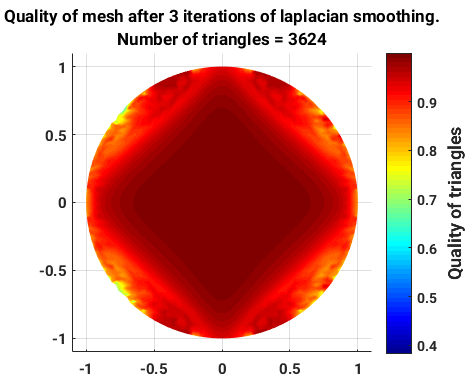}
\caption{\bl{Quality of the mesh after 3 iterations of Laplace smoothing applied to the original un-smoothed deformed mesh from Figure \ref{fig:1c}}}
\label{fig:fig_laplacian1}  
\end{figure} 
Figure \ref{fig:1b} shows a surface plot of the per-vertex quality measure on the undeformed domain. Figure \ref{fig:1d} shows the per-vertex quality in the deformed domain obtained from RBF interpolation, and the corresponding mesh is shown in Figure \ref{fig:1c}. As predicted, the deformed mesh shows mild distortions corresponding to the singularities in the conformal map. Our smoothing algorithm terminates after three iterations with little improvement in per-vertex quality (Figures \ref{fig:1e} and \ref{fig:1f}). In fact, Figure \ref{fig:figuretype1c} shows that the per-vertex quality \emph{decreases} after just three iterations, which is indeed the cause for termination. We posit that this lack of improvement in quality is due to the fact that the deformed mesh is of relatively high quality to begin with. It is likely that a purely local algorithm would be able to achieve higher per-vertex quality than our quasi-local algorithm. \bl{Figure \ref{fig:fig_laplacian1} shows the result of applying three iterations of Laplace smoothing to the deformed, unsmoothed mesh from Figure \ref{fig:1c}. Laplace smoothing, being a purely local technique, improves the quality of elements faster than the RBF-based smoothing. }

\subsection{Deforming an annulus into a square with an airfoil cavity}
\begin{figure}[h!]
\centering
\subfigure[Initial tessellated undeformed domain]{%
 \includegraphics[width=.43\textwidth]{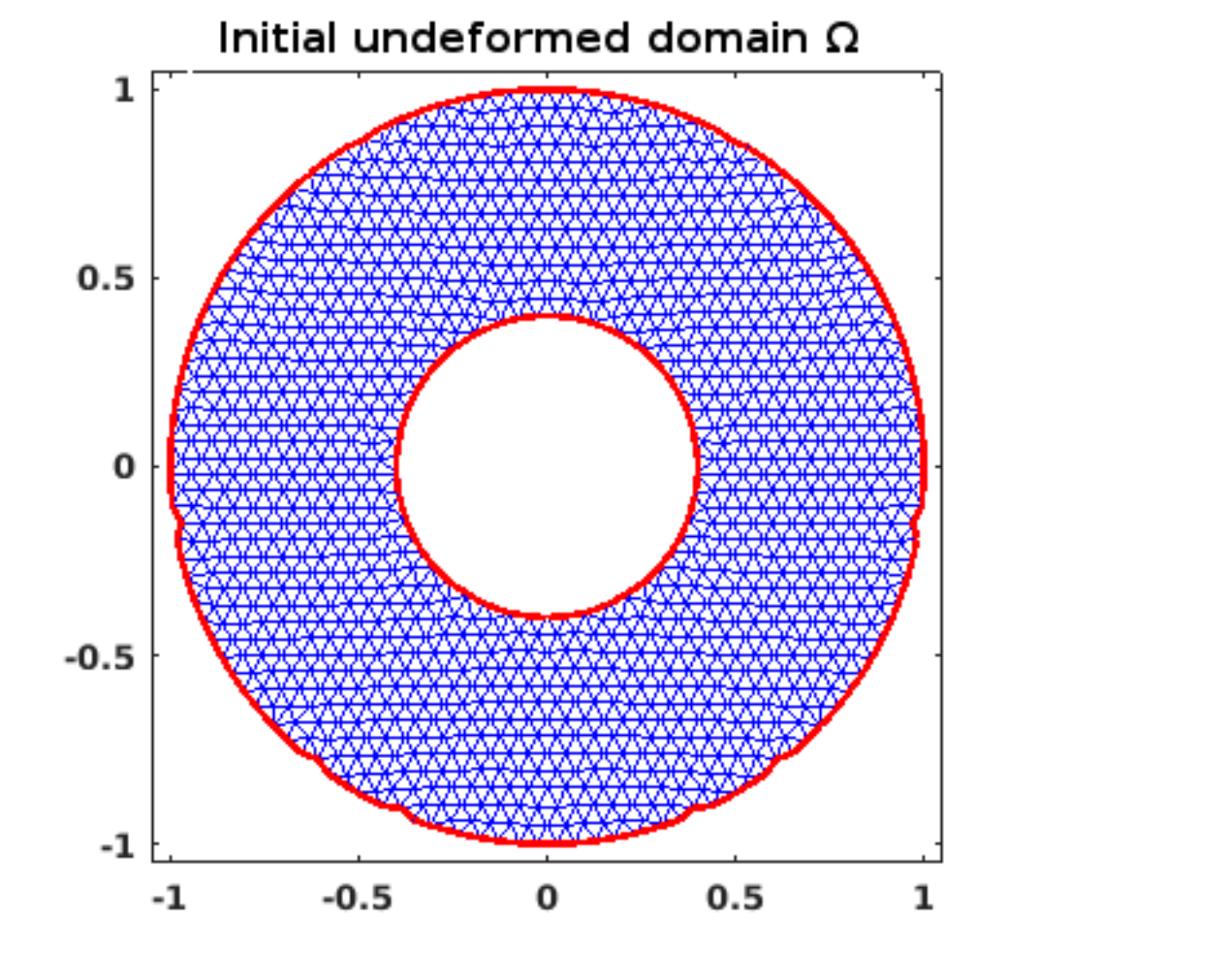}
\label{fig:2a}}
\subfigure[Quality of undeformed initial mesh]{%
 \includegraphics[width=.43\textwidth]{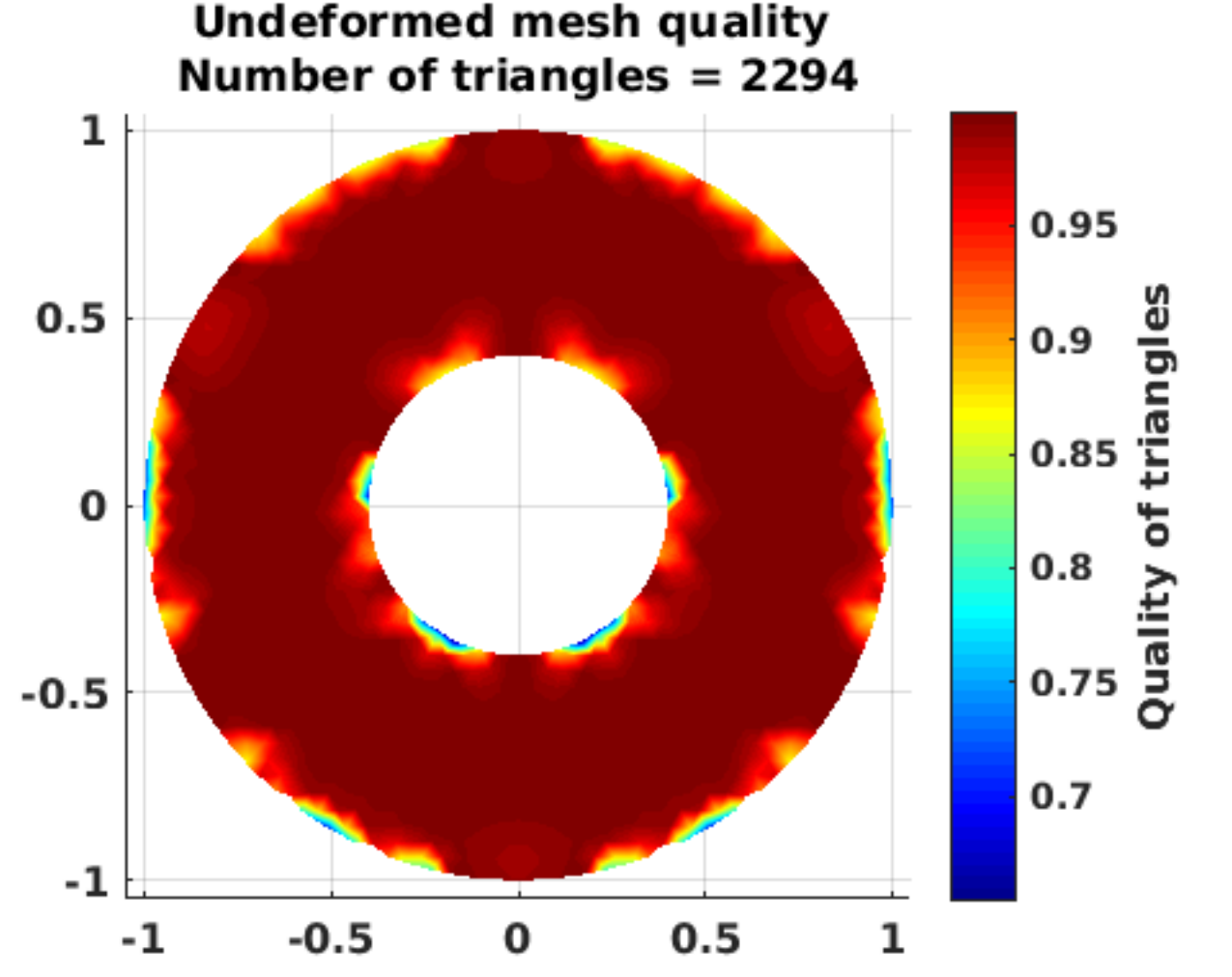}
\label{fig:2b}}

\subfigure[Deformed mesh before smoothing]{
\includegraphics[width=.43\textwidth]{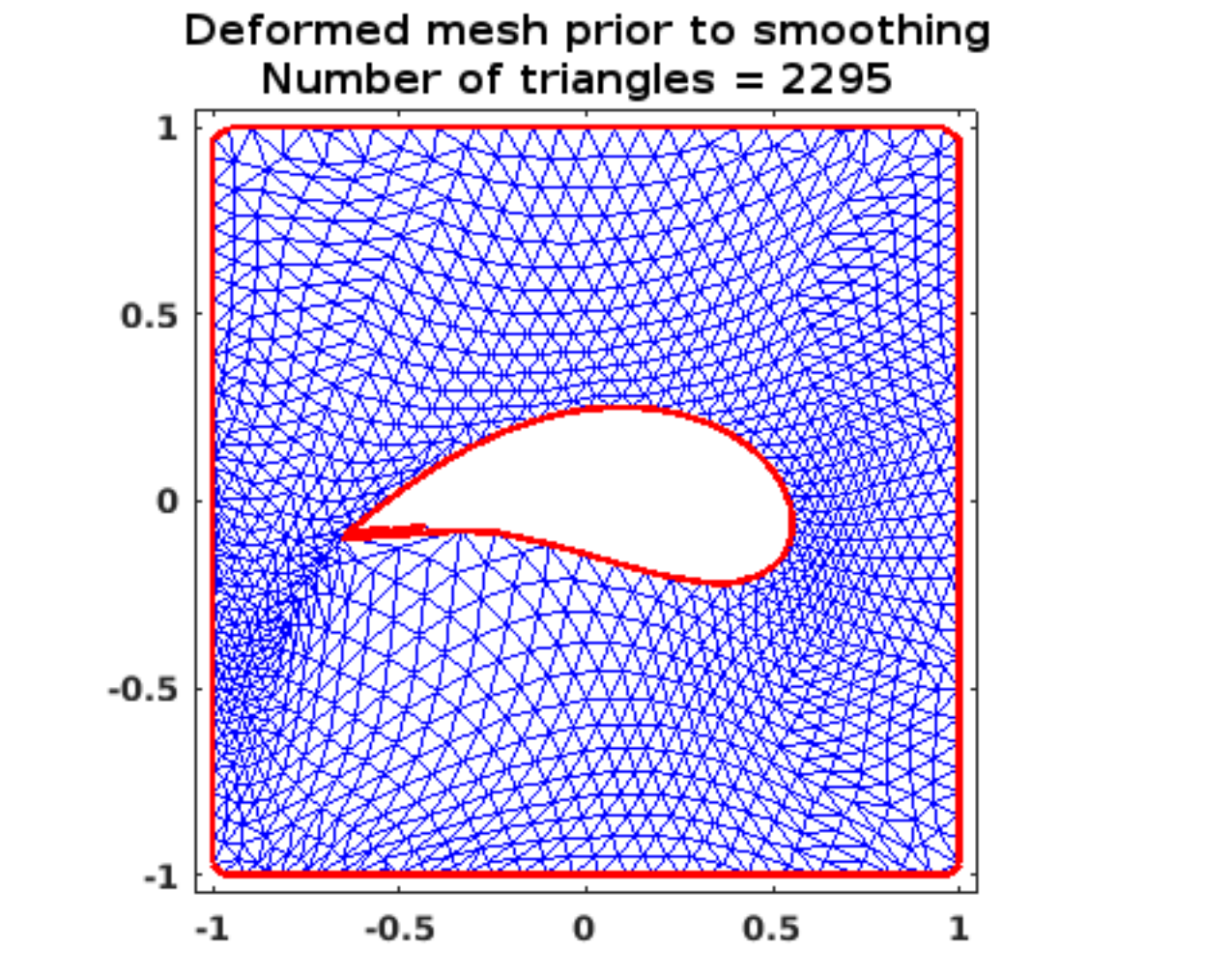}
\label{fig:2c}}
\subfigure[Quality of deformed mesh before smoothing]{%
 \includegraphics[width=.43\textwidth]{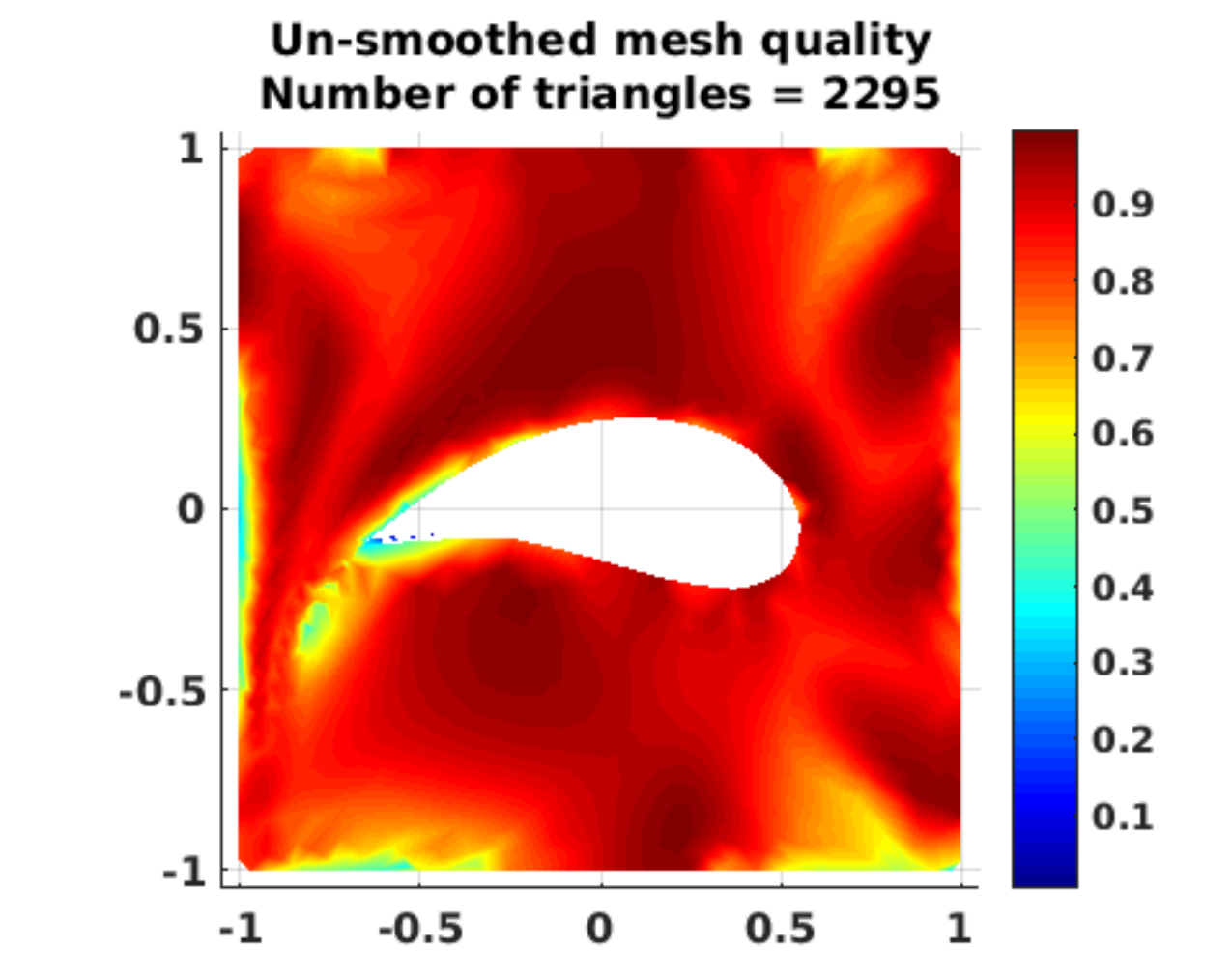}
\label{fig:2d}}

\subfigure[Deformed mesh after smoothing]{%
 \includegraphics[width=.43\textwidth]{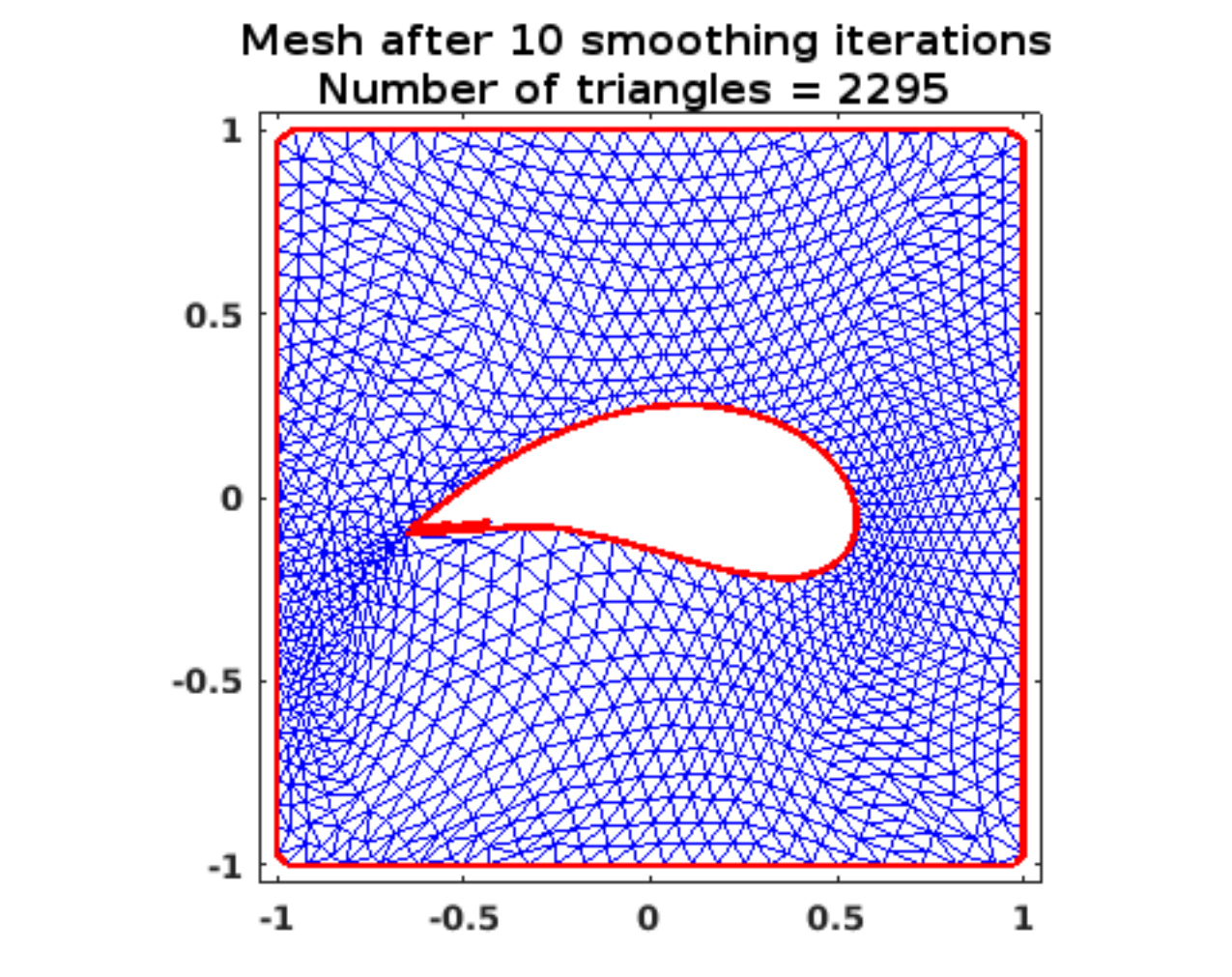}
\label{fig:2e}}
\subfigure[Quality of deformed mesh after smoothing]{%
 \includegraphics[width=.43\textwidth]{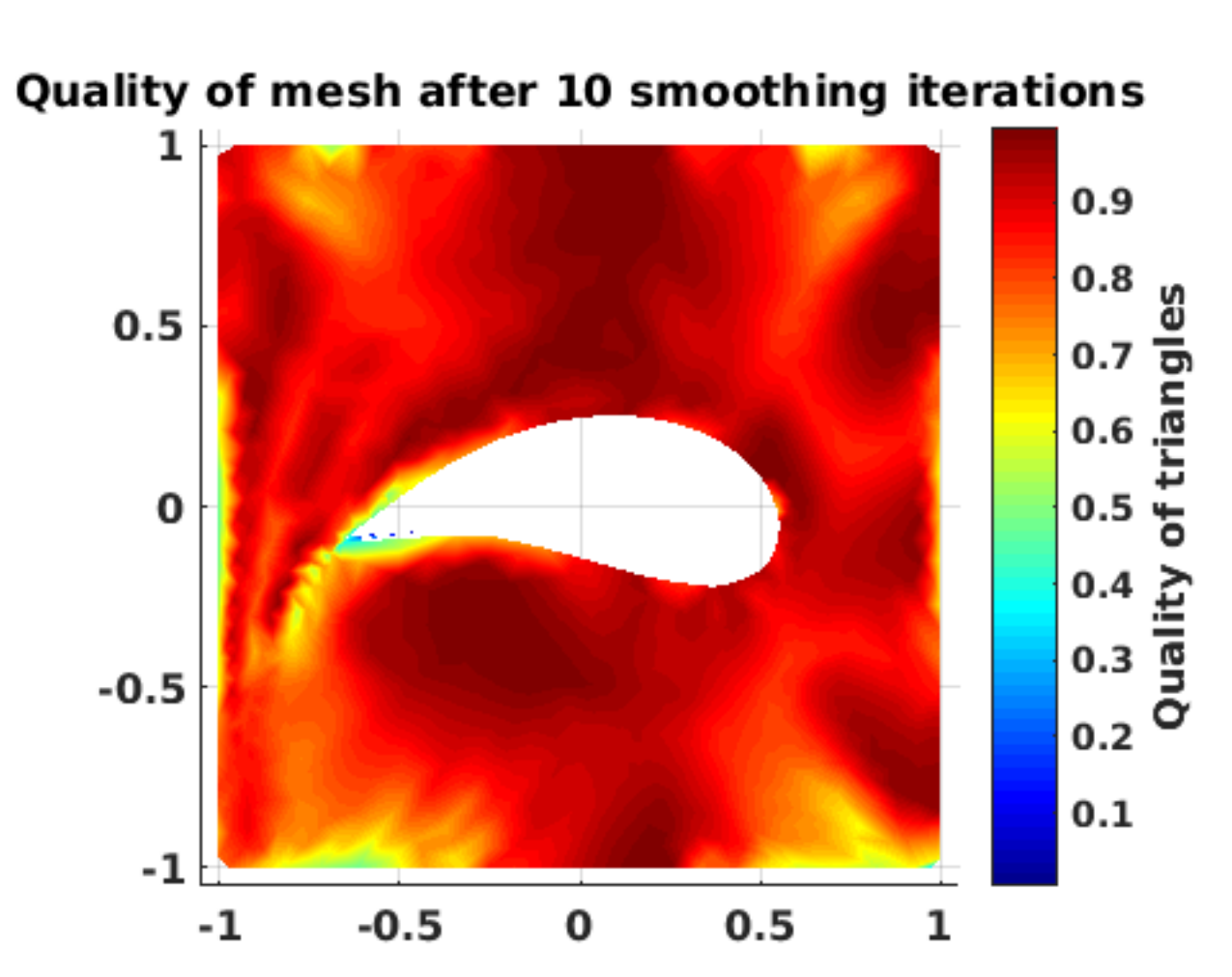}
\label{fig:2f}}
\caption{Deforming an annulus into a square with an airfoil cavity. }
\label{fig:figuretype2b}  
\end{figure} 
We now consider a more complicated test case. In this test, we deform a circular annulus into a square containing an airfoil. Essentially, this test transforms both inner and outer boundaries from $C^{\infty}$ smoothness to $C^1$ smoothness. This test case shows the ability of our method to naturally handle embedded boundaries. The deformation function on the cavity boundary is the standard Joukowsky conformal map\cite{tsien1943symmetrical} from circle to airfoil, while the deformation function on the outer boundary is another conformal map from a circle to a square. The results of this experiment are shown in Figure \ref{fig:figuretype2b}.
\begin{figure}[h!]
\centering
\includegraphics[width=.43\textwidth]{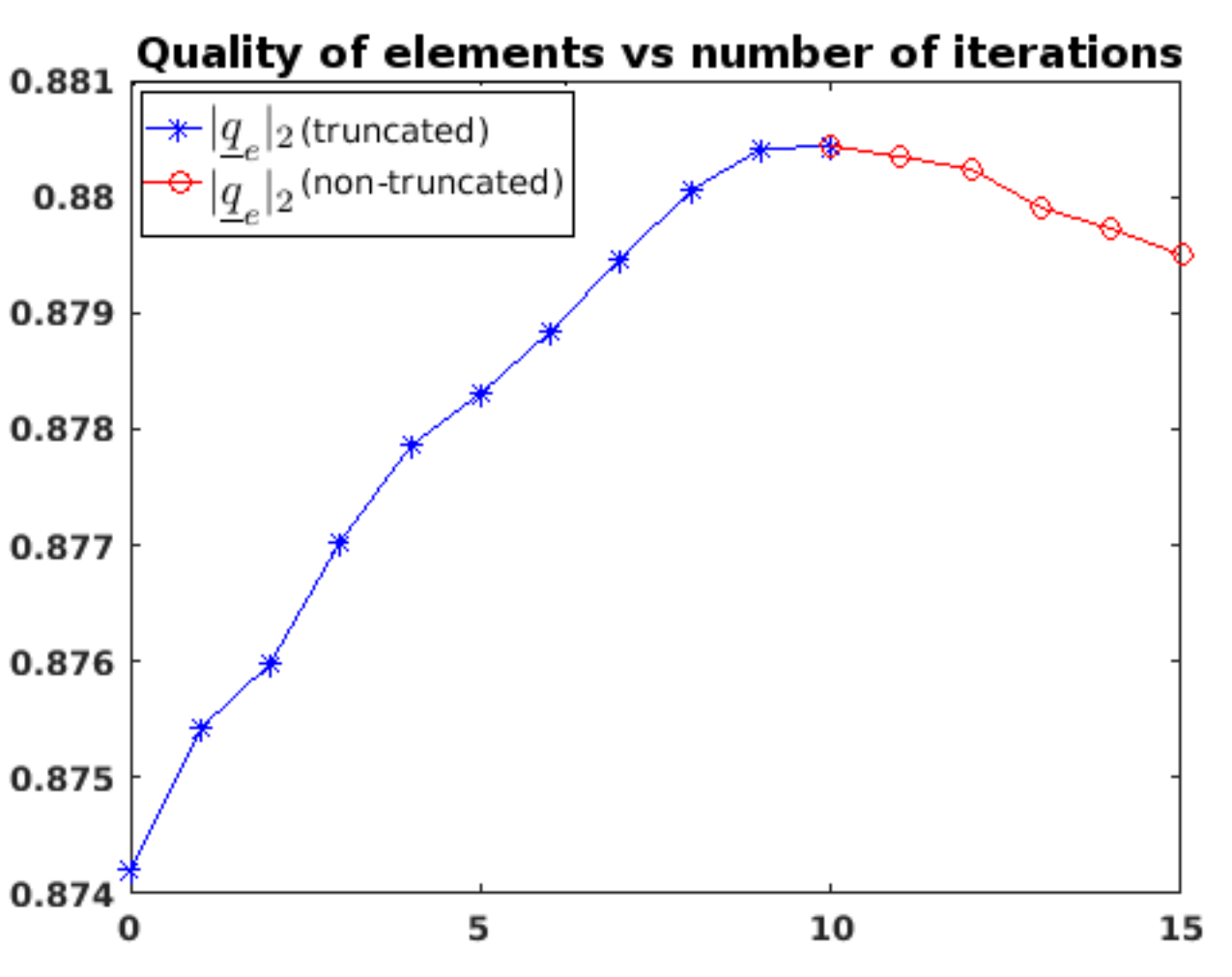}
\caption{Quality as a function of number of iterations when deforming an annulus into a square with an inner airfoil. The blue line corresponds to the use of a stopping criterion. The red line shows the quality when no stopping criterion is used.}
\label{fig:figuretype2c}  
\end{figure} 
\begin{figure}[h!]
 \centering
 \includegraphics[width=.50\textwidth]{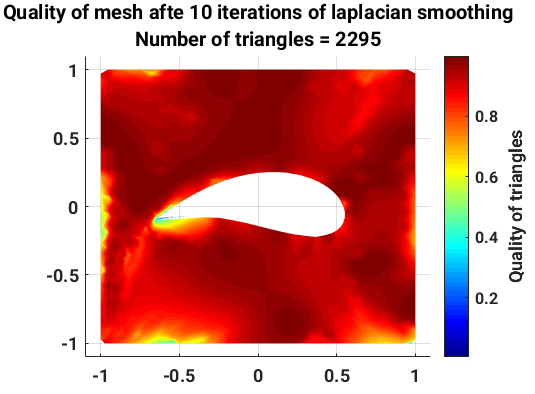}
 \caption{\bl{Quality of the mesh after 10 iterations of Laplace smoothing applied to the original un-smoothed deformed mesh from Figure \ref{fig:2c}}}
 \label{fig:fig_laplacian2}  
\end{figure}
Figure \ref{fig:2b} shows the per-vertex quality surface plot of the undeformed mesh, while Figure \ref{fig:2d} shows the same plot for the deformed mesh before smoothing. The singularity in the Joukowsky map manifests as poor quality elements near the trailing edge of the airfoil. Our smoothing algorithm now runs for ten iterations before terminating; the resulting mesh is shown in Figure \ref{fig:2e}, and its quality is shown in Figure \ref{fig:2f}. The benefits of smoothing are apparent: the mesh distortions near the trailing edge have been reduced without adversely affecting the higher-quality regions. The mesh distortion near the left edge of the square has also been reduced. A study of the mesh quality as a function of the number of iterations is shown in Figure \ref{fig:figuretype2c}. While it is impossible to capture overall mesh quality with a single number, it is easy to see the improvement of quality due to the smoothing algorithm and the benefit of the termination criterion. \bl{Figure \ref{fig:fig_laplacian2} shows the result of applying ten iterations of Laplace smoothing to the deformed, un-smoothed mesh from Figure \ref{fig:2c}. Once again, Laplace smoothing performs the smoothing faster due to its local nature.}

\subsection{Deforming a cube to a sphere}
\begin{figure}[h]
\centering
\subfigure[Undeformed Domain]{
\includegraphics[width=.43\textwidth]{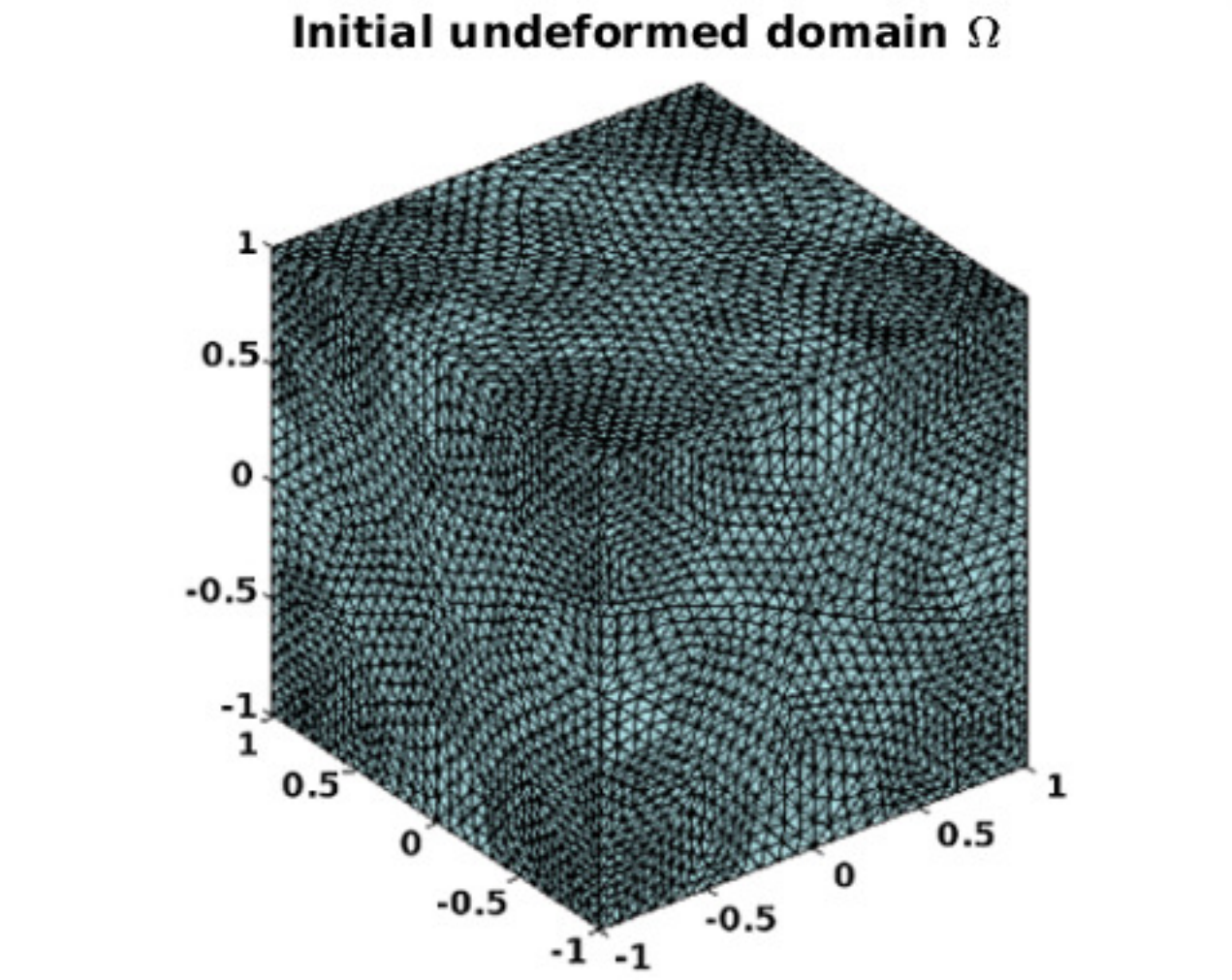}\label{fig:4a}
}
\subfigure[Exterior view of element quality]{
\includegraphics[width=.43\textwidth]{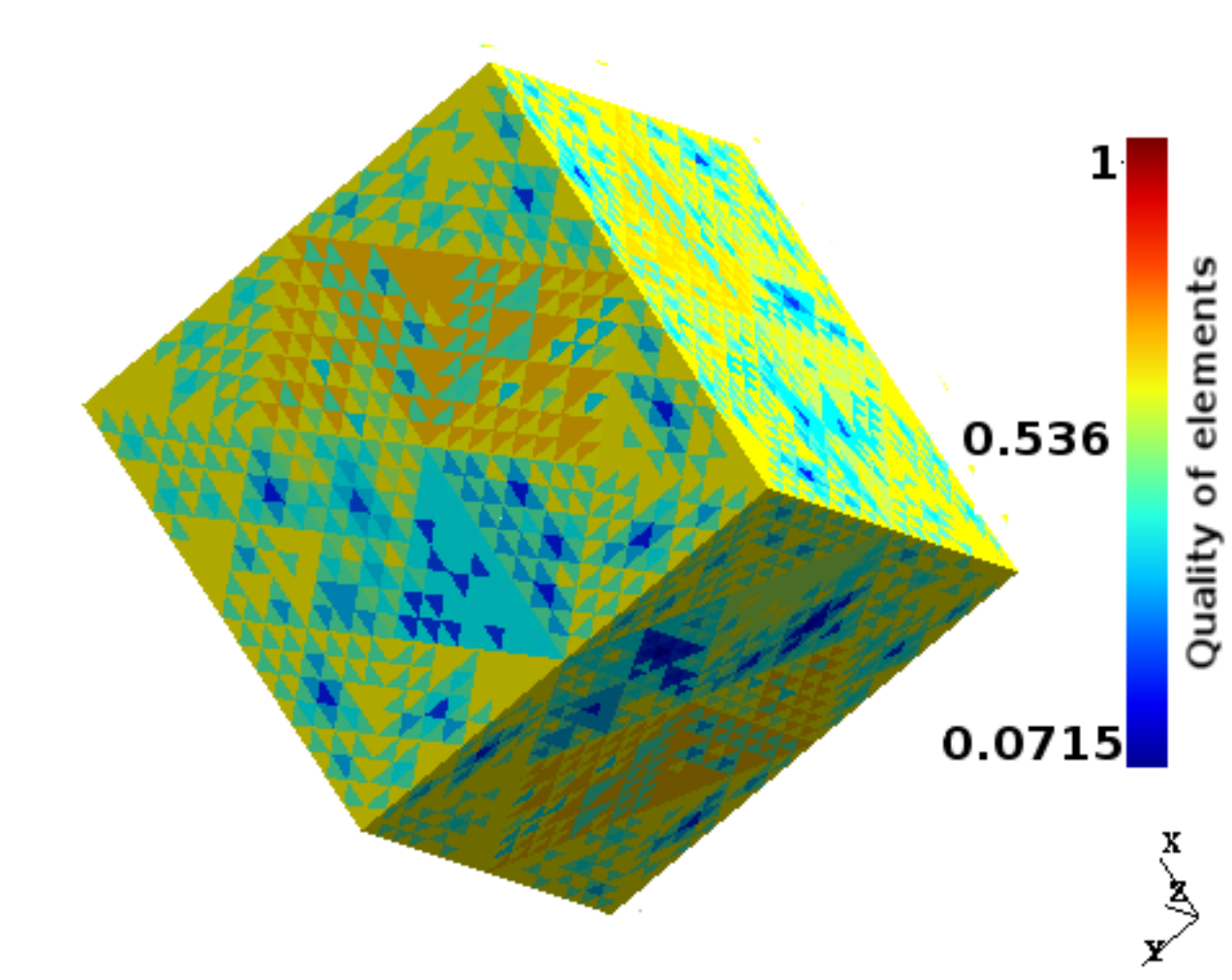}\label{fig:4c}
}
\subfigure[Interior view of element quality]{
\includegraphics[width=.43\textwidth]{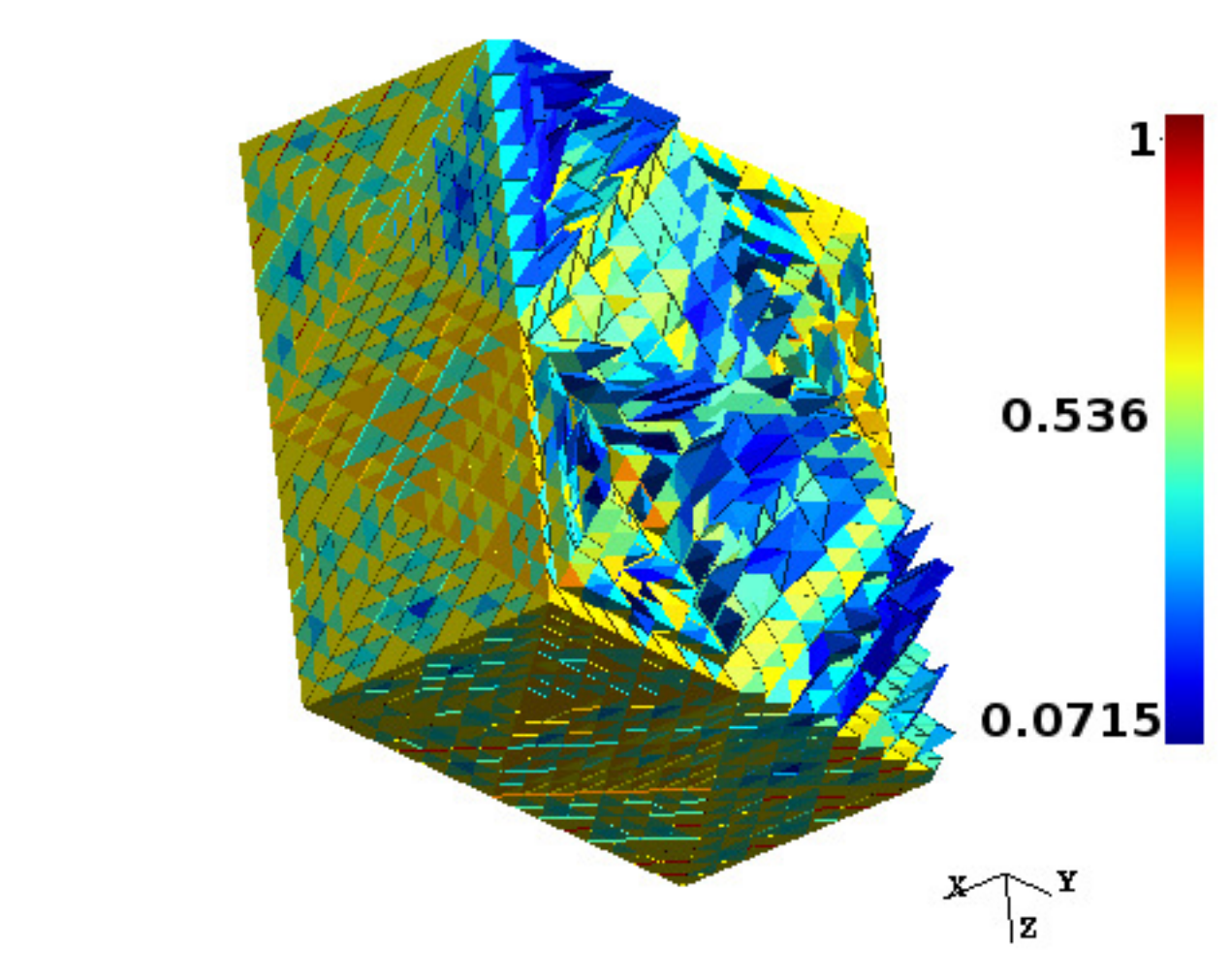}\label{fig:4d}
}
\caption{Tetrahedral cube mesh and element quality.}
\label{fig:4all}
\end{figure}
\begin{figure}[h]
\centering
\subfigure[Deformed domain]{
\includegraphics[width=.43\textwidth]{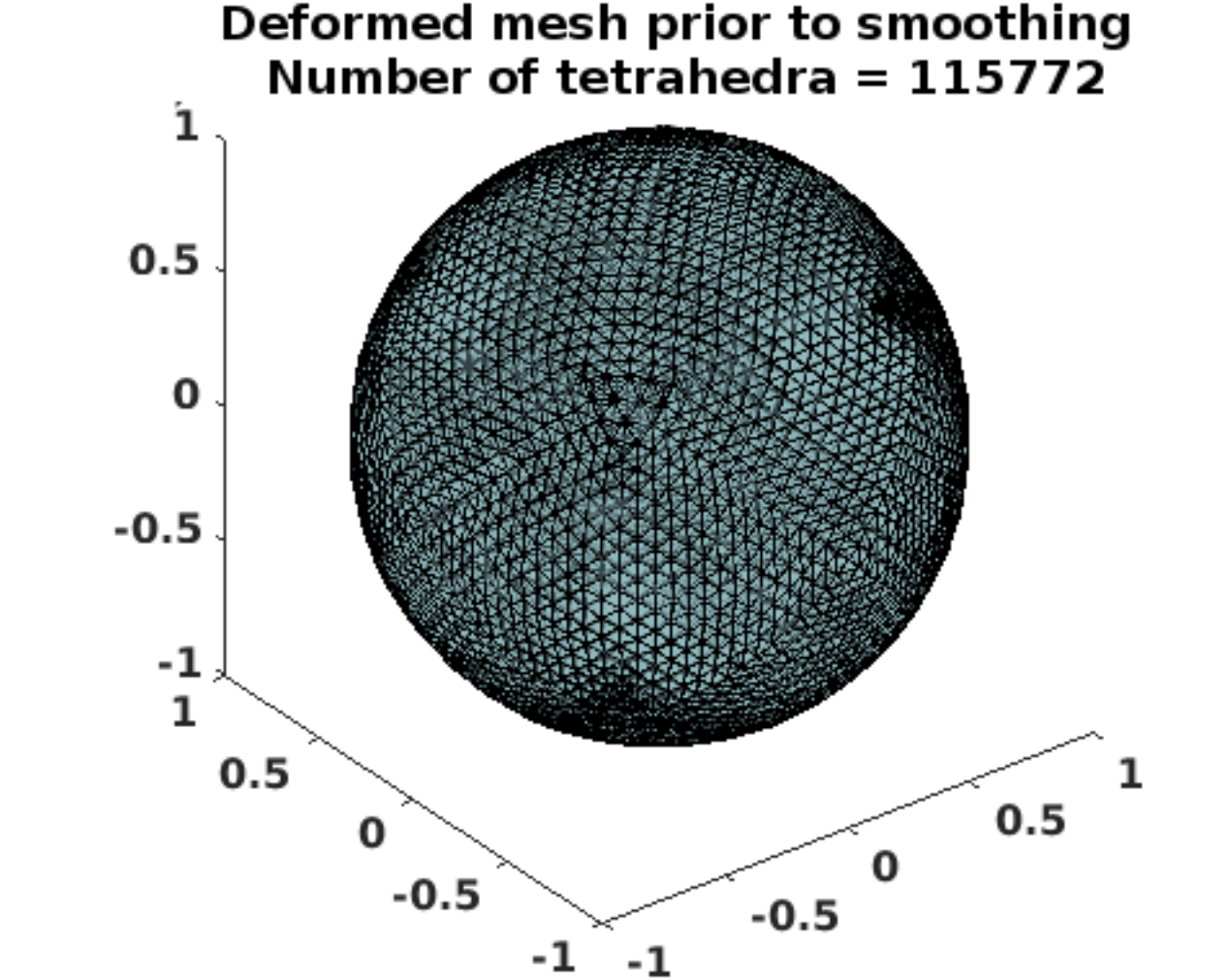}\label{fig:5a}
}
\subfigure[Exterior view of element quality]{%
 \includegraphics[width=.43\textwidth]{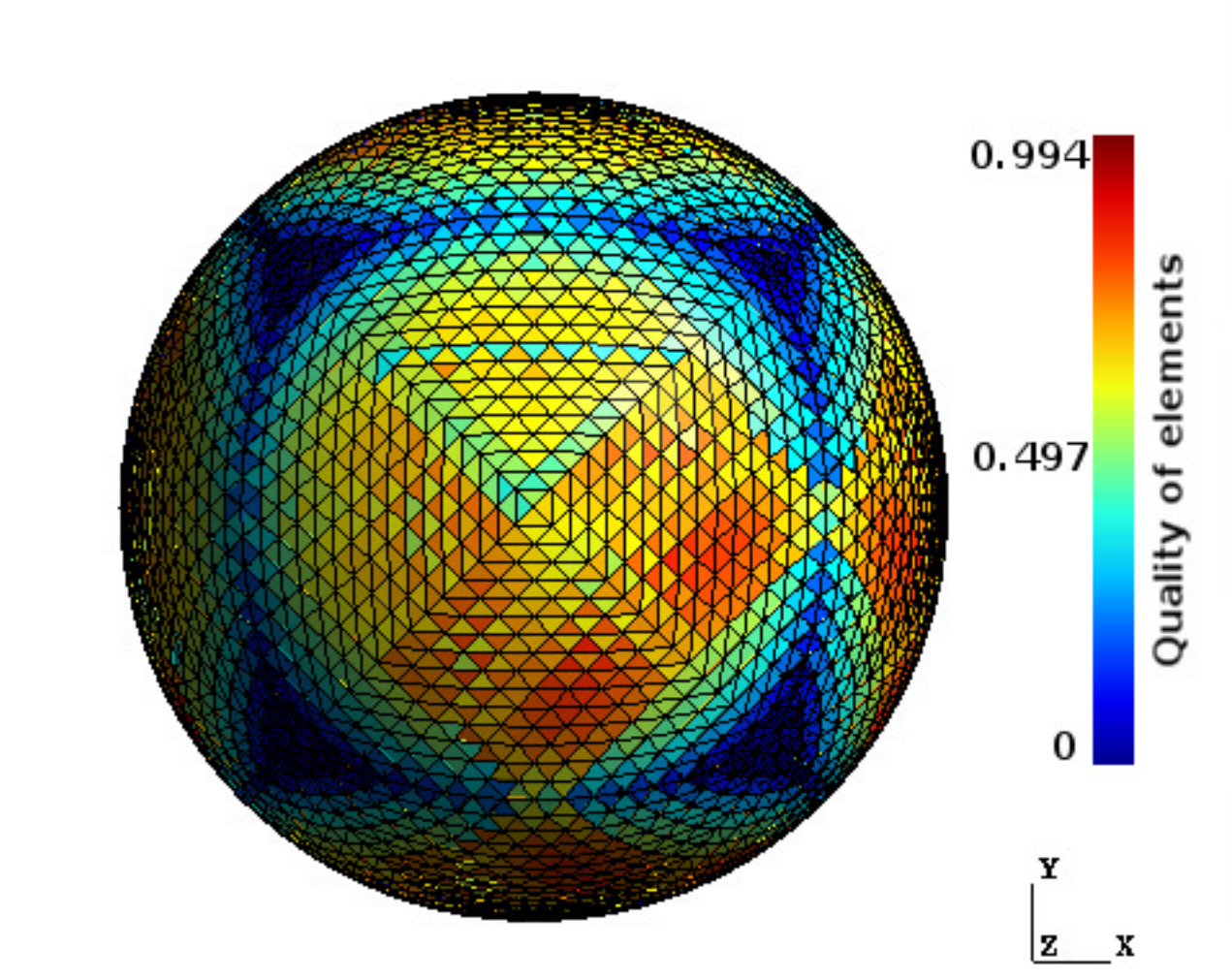}
\label{fig:5b}}

\subfigure[Interior view of element quality]{%
 \includegraphics[width=.43\textwidth]{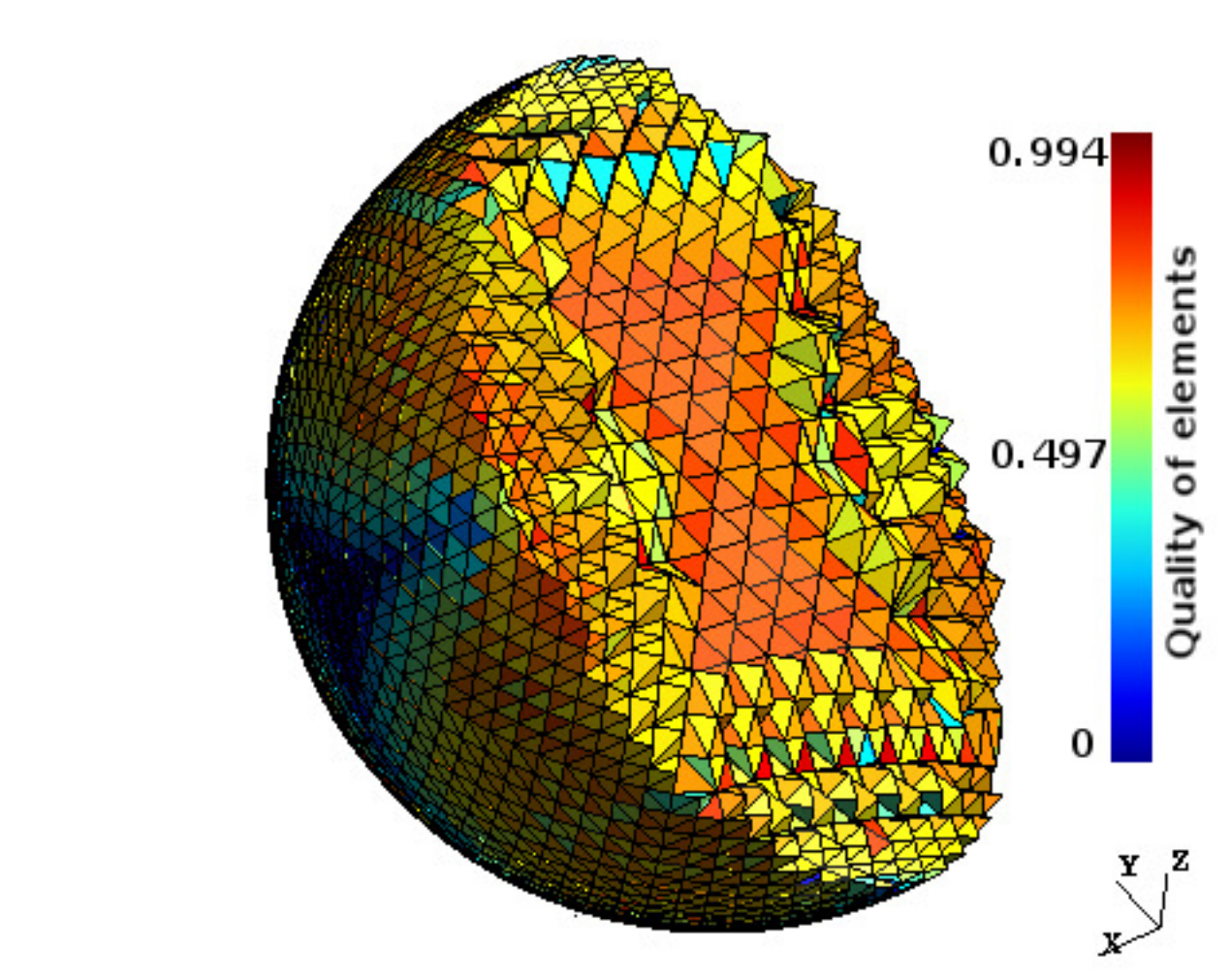}
\label{fig:5c}}
\subfigure[Quality as a function of number of iterations when deforming a unit cube to a sphere. The blue line corresponds to the use of a stopping criterion. The red line shows the quality when no stopping criterion is used.]{%
 \includegraphics[width=.43\textwidth]{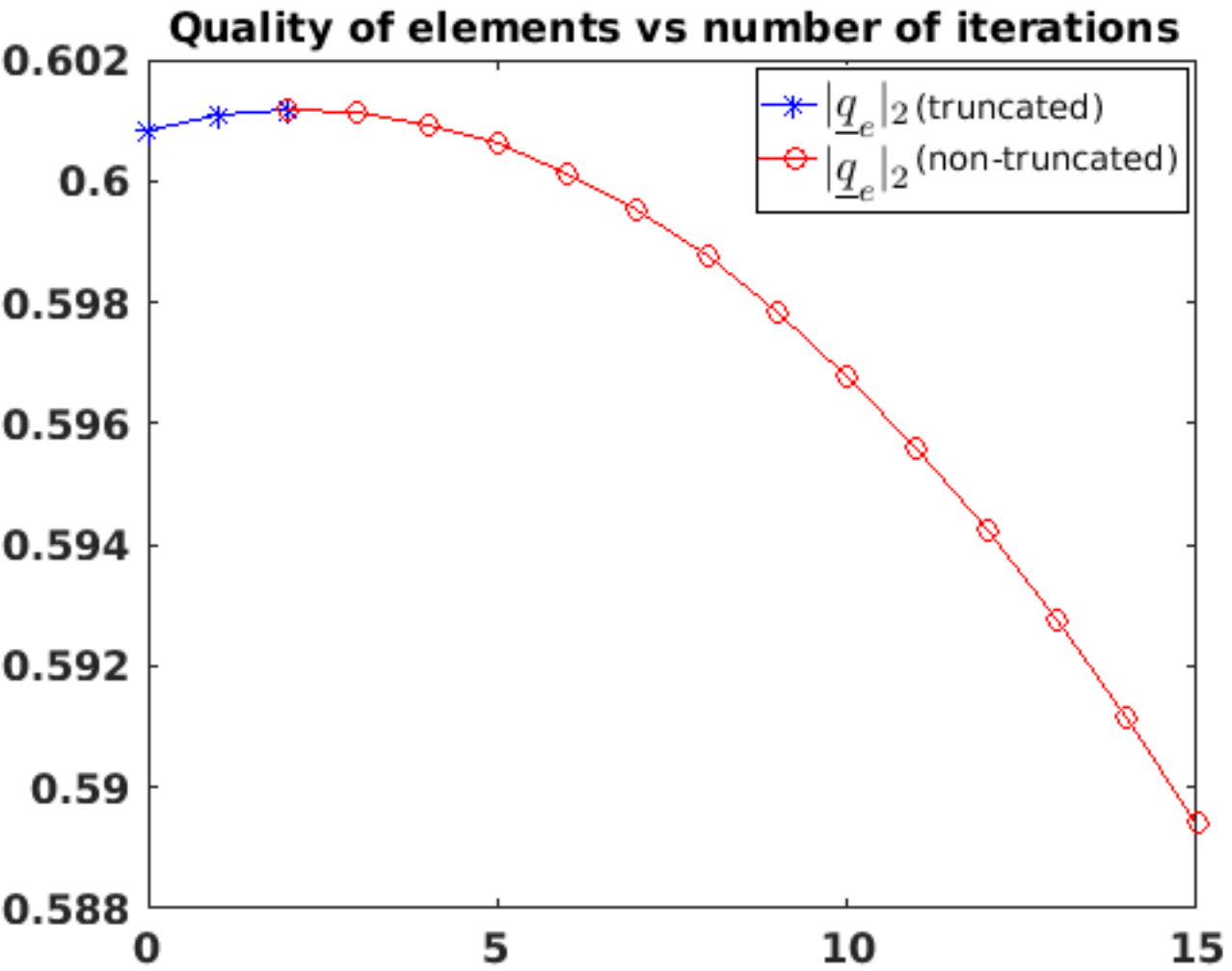}
\label{fig:5d}}
\caption{Curvilinear mesh obtained by deforming the cube to the sphere.}
\label{fig:5all}
\end{figure}

One of the main advantages of RBF interpolation is that it requires no modifications for interpolating data scattered on submanifolds $\mathbb{M} \subset \mathbb{R}^3$~\cite{FuselierWright2012}. Thus, our algorithm requires no modification in $\mathbb{R}^3$ beyond recovering the deformation map for the third coordinate. We work completely in Cartesian coordinates and do not employ special node sets on our undeformed domain boundary. We now consider the 3D analogue of the square-to-disk test by deforming the unit cube to a sphere. Once again, as in the 2D test cases, we employ a straightforward conformal map from the cube to the sphere. The undeformed domain mesh is shown in Figure \ref{fig:4a}, and Figures \ref{fig:4c} and \ref{fig:4d} show the element quality in that mesh. It is easy to see from Figures \ref{fig:4c} and \ref{fig:4d} that the mesh is mostly comprised of low quality elements on the boundary, and very low quality elements in the interior. We now apply the RBF interpolation and smoothing procedure to the mesh to obtain a mesh within a sphere. The resulting mesh and element quality after smoothing are shown in Figure \ref{fig:5all}.
Figure \ref{fig:5a} shows the curvilinear mesh obtained with the sphere. An exterior view of element quality (Figure \ref{fig:5b}) shows that the mesh quality is lowest at spatial locations corresponding to singularities in the conformal map. Interestingly, both Figures \ref{fig:5b} and \ref{fig:5c} show that the overall element quality in the deformed mesh is higher than in the undeformed mesh, especially in the interior. This illustrates a strength of the RBF technique in generating curvilinear meshes. It should be noted that the improvement in quality by smoothing does not invalidate any elements in the mesh. This can be easily verified by observing that no element in the deformed mesh has a negative jacobian. Finally, the last plot shows that our smoothing procedure terminates after just two iterations as further improvement is not possible. Indeed, this test case indicates that the quality of the undeformed mesh proves crucial in dictating the quality of the curvilinear smoothed mesh. However, this test case still illustrates that our technique is viable in 3D without any real changes to the algorithm. \bl{We also tested our technique using a high-quality mesh on the undeformed domain. In this case, we note that our smoothing algorithm did not significantly improve the quality of the mesh on the deformed domain (results not shown).}

\section{Discussion}
\label{sec:conclusion}

The main contribution of this article is a framework for generating 2D and 3D curvilinear meshes using RBF interpolation on the domain boundary, and a quasi-local iterative algorithm for smoothing those meshes by modifying RBF shape parameters. Interestingly, the technique allows mesh generation in the interior of the domain using an approximation to the deformation map built solely on the domain boundary. Despite the maps not being harmonic functions, this technique appears to produce meshes that either preserve or improve the quality of the undeformed mesh. Our results indicate that smoothing can be beneficial, especially with meshes produced from singular deformation maps such as the Joukowsky transform. Further, our algorithm is directly applicable to both 2D and 3D mesh generation in Cartesian coordinates due to the ability of RBF interpolants to handle scattered node sets on submanifolds of $\mathbb{R}^d$.

Despite its quasi-local nature, our smoothing algorithm is still global due to the use of global interpolants on the boundary. This likely limits the ability of our algorithm to handle isolated low-quality regions without adversely affecting high-quality regions. A natural approach to overcome this will be to use RBF-based Partition of Unity (RBF-PU) or RBF-based Finite Difference (RBF-FD) methods to approximate the boundary deformation map. These methods retain the advantages of global RBF interpolants-- the ability to handle scattered data on submanifolds of $\mathbb{R}^d$, high-order convergence rates for smooth functions-- but have lower costs and are more localized. Both these methods have the potential to bring down the preprocessing costs from $O(N^{1.5})$ and $O(N^2)$ to $O(N)$. Further, per-iteration evaluation costs can also be decreased to $O(N)$ from $O(N^{\frac{3}{2}})$ and $O(N^{\frac{5}{3}})$. This is an area of future research.

We remark that our algorithm is likely applicable in many scenarios beyond generating curvilinear meshes. For instance, we envision that our algorithm may be useful in remeshing particle meshes in Lagrangian methods, or in generating node sets for the meshfree or meshed solution of PDEs on time-varying domains. Finally, an open area of research is to understand how RBF interpolation on the boundary is able to recover a non-harmonic conformal map in the interior of the domain.

\begin{acknowledgements}
VZ was supported by NSF OCI-1148291 and NSF IIS-1212806. VS was supported by NSF DMS-1521748. SPS was supported in part by the NIH/NIGMS Center for Integrative Biomedical Computing grant 2P41 RR0112553-12 and a grant from the ExxonMobil corporation. RMK was supported in part by DMS-1521748 and W911NF-15-1-0222.
\end{acknowledgements}

\bibliographystyle{spmpsci}      % mathematics and physical sciences
\bibliography{Article}
\end{document}